\documentclass[12pt,letterpaper]{amsart}
\usepackage{amsmath,amssymb,amsfonts,amsthm}
\usepackage{mathrsfs}
\usepackage[all]{xy}
\usepackage{stmaryrd}
\usepackage{mathtools}
\usepackage{indentfirst}
\usepackage{verbatim}
\usepackage{hyperref}
\usepackage{color}
\usepackage[usenames,dvipsnames]{xcolor}
\usepackage{graphicx}
\graphicspath{ {images/} }
\usepackage[normalem]{ulem}
\usepackage{tikz}
\usepackage{nicefrac}
\usepackage{xfrac}
\usepackage[shortlabels]{enumitem}

\title{Log canonical foliation singularities on surfaces}
\author{Yen-An Chen}
\subjclass[2010]{Primary 32S65, Secondary 32M25, 14B05, 14J99.}
\thanks{The author was partially supported by NSF research grants no: DMS-1801851, DMS-1840190 and by a grant from the Simons Foundation; Award Number: 256202. }
\address{Department of Mathematics, University of Utah, Salt Lake City, UT 84112, USA}
\email{yachen@math.utah.edu}

\newtheorem*{claim}{Claim}
\newtheorem{thm}{Theorem}[section]
\newtheorem{prop}[thm]{Proposition}
\newtheorem{cor}[thm]{Corollary}
\newtheorem{lem}[thm]{Lemma}

\theoremstyle{definition}
\newtheorem{defn}[thm]{Definition}
\newtheorem*{acks}{Acknowledgements}

\theoremstyle{remark}
\newtheorem{rmk}[thm]{Remark}
\newtheorem*{pf}{Proof}

\newcommand\cE{{\mathcal{E}}}

\newcommand\cO{{\mathcal{O}}}

\newcommand\bC{{\mathbb C}}

\newcommand\bN{{\mathbb N}}
\newcommand\bP{{\mathbb P}}
\newcommand\bQ{{\mathbb Q}}
\newcommand\bR{{\mathbb R}}

\newcommand\bZ{{\mathbb Z}}

\newcommand\sF{{\mathscr{F}}}
\newcommand\sG{{\mathscr{G}}}
\newcommand\sH{{\mathscr{H}}}

\newcommand\rw{\rightarrow}
\newcommand\wt[1]{\widetilde{#1}}

\newcommand\ve{\varepsilon}

\newcommand{\D}{\Delta}
\newcommand{\T}{\Theta}
\newcommand{\G}{\Gamma}
\newcommand{\bs}{\backslash}

\begin{document}

\maketitle

\begin{abstract}
We give a classification of the dual graphs of the exceptional divisors on the minimal resolutions of log canonical foliation singularities on surfaces. 
For an application, we show the set of foliated minimal log discrepancies for foliated surface triples satisfies the ascending chain condition and a Grauert-Riemenschneider type vanishing theorem for foliated surfaces with good log canonical foliation singularities. 
\end{abstract}

\section*{Introduction}

Singularities play an important role in many areas of algebraic geometry. 
For instance, singular varieties naturally appear in the minimal model program and in the study of moduli spaces, where it often happens that smooth objects degenerate to singular ones. 
However, many classes of singularities that naturally occur are mild and can sometimes be understood in detail.
One of the natural classes of singularities is (log) canonical singularities. 
These singularities have been extensively studied, especially in the case of surfaces.
In particular, we have a full classification of log canonical surface singularities. 
(For a reference, see \cite{watanabe1980plurigenera}, \cite[Section 4.1]{kollar1998birational}, or \cite[Appendix]{hara1998classification}.) 

When studying foliations, one expects that similar results may hold. 
In this direction, \mbox{McQuillan} introduces in \cite{mcquillan2008canonical} a notion of log canonical foliation singularities, which is a natural generalization of log canonical singularities. 
Moreover, he also gives a classification for the canonical foliation singularities on surfaces. 
(See \cite[Corollary I.2.2 and Fact I.2.4]{mcquillan2008canonical}.) 

In this paper, we first give a full list of all possibilities of the dual graphs of the exceptional divisors on the minimal resolutions of log canonical foliation singularities on surfaces. 
(See also Theorem~\ref{main}.) 
This is achieved with the help of the work of Brunella, Camacho, and Sad on the indices and a theorem on separatrices. 

\begin{thm}
Let $(X,\sF,p)$ be a germ of a foliated surface. 
Assume that $p$ is a log canonical singularity of $\sF$. 
Let $\pi : (Y,\sG) \rightarrow (X,\sF)$ be the minimal resolution (see subsection~\ref{min_rsln_fol}) for $(X,\sF,p)$ with exceptional divisors $E = \cup E_i$. 
Then $E$ belongs to one of the following types: 
\begin{enumerate}
\item A $\sG$-chain. 
\item A chain of three invariant curves $E_1\cup E_2\cup E_3$ where $E_1$ and $E_3$ are $(-1)$-$\sG$-curves with self-intersection $-2$ and $E_2$ is a bad tail.
\item A chain of $(-2)$-$\sG$-curves. 
\item A dihedral singularity. 
More precisely, two $(-1)$-$\sG$-curves with self-intersection $-2$ joined by a bad tail which itself connects to a chain of $(-2)$-$\sG$-curves. 
\item An elliptic Gorenstein leaf. 
\item A chain $E = \bigcup_{i=1}^rE_i$ with exactly one non-invariant curve $E_\ell$ with $1\leq \ell\leq r$. Moreover, $E_\ell$ has tangency order zero and $\bigcup_{i=1}^{\ell-1}E_i$ and $\bigcup_{i=\ell+1}^{r}E_i$ are $\sG$-chains. 
\item The dual graph is star-shaped with a non-invariant center $[E_0]$. Moreover, $E_0$ has tangency order zero, all branches are $\sG$-chains, and all first curves of $\sG$-chains have intersection number one with $E_0$. 
\end{enumerate}
Note that type (1) is terminal, and types (1) - (5) are canonical. 
\end{thm}

Inspired by the work in \cite{alexeev1993two} on the ascending chain condition (ACC) for the set of the minimal log discrepancies of surface singularities, we show the set of foliated minimal log discrepancies satisfies the ACC. 
(See also Definition~\ref{foliated_discrepancy} and Theorem~\ref{mld_acc}.)
\begin{thm}
For any set $B$ satisfying the descending chain condition, the set 
\[\textnormal{MLD}(2,B) := \{\textnormal{mld}_x(\sF,\Delta) \vert \, \mbox{$(X,\sF,\Delta)$ is a foliated triple with $x\in X$ and $\Delta\in B$}\}\] 
satisfies the ascending chain condition (ACC). 
\end{thm}

Finally, we prove a Grauert-Riemannschneider type vanishing theorem for foliated surfaces with \emph{good} log canonical foliation singularities by using the method in \cite[Theorem 10.4]{kollar2013singularities}. 
(See Definition~\ref{good_log_canonical} and Theorem~\ref{van_thm}.)
This is a generalization of \cite[Theorem 6.1]{hacon2021birational} in which $(X,\sF)$ is assumed to have only canonical foliation singularities.

\begin{thm}
Let $f : (Y,\sG) \rw (X,\sF)$ be a proper birational morphism where $(X,\sF)$ is a foliated surface with good log canonical foliation singularities and $(Y,\sG)$ is a foliated surface with only reduced singularities. 
Then $R^if_*\cO_Y(K_\sG)=0$ for $i>0$. 
\end{thm}

\begin{acks}
The author would like to thank Christopher D. Hacon for his insightful suggestions and encouragements. 
\end{acks}

\section{Preliminaries}
In this paper, we always work over $\bC$. 
By a surface, we mean a normal algebraic space of dimension two. 
In this section, we recall several definitions and results which will be used later. 

\subsection{Foliations on surfaces}
A \emph{foliation} $\sF$ on a surface $X$ is a rank 1 saturated subsheaf $\sF$ of the tangent sheaf $T_{X}$ of $X$. 
So we have the following short exact sequence: 
\[0\rw \sF\rw T_X \rw T_X/\sF \rw 0\] 
with $T_X/\sF$ torsion-free. 

The point $p$ on $X$ is called a \emph{singular} point of the foliation if either a singular point of $X$ or a point at which the quotient $T_X/\sF$ is \emph{not} locally free. 
Since $X$ is normal, the foliation singularities are isolated. 
\begin{defn}
A \emph{foliated surface} is a pair $(X,\sF)$ consisting of a surface $X$ and a foliation $\sF$ on $X$. 
A \emph{foliated triple} is a triple $(X,\sF,\D)$ consisting of a foliated surface $(X,\sF)$ and an $\bR$-divisor $\D=\sum a_iD_i$. 
\end{defn}
Notice that $T_X\cong \textnormal{Hom}_{\cO_X}(\Omega_X,\cO_X)$ is reflexive. 
Thus, $\sF$ is also reflexive and therefore we can define the \emph{canonical divisor} $K_\sF$ of the foliation as a Weil divisor on $X$ with $\cO_X(-K_\sF)\cong \sF$. 

\begin{defn}
Let $(X,\sF)$ be a foliated surface. 
Given any birational map between normal surfaces $f : Y\dashrightarrow X$ and a foliation $\sF$ on $X$, then we define the \emph{pullback foliation} $f^*\sF$ as follows: 

Let $U$ be an open subset such that $V:= f^{-1}(U) \rightarrow U$ is an isomorphism. 
Note that $\sF\vert_U \subset T_U \cong T_V$. 
By \cite[Exercise II.5.15]{hartshorne1977algebraic}, we have a coherent subsheaf $\sG$ of $T_Y$ such that $\sG\vert_V = \sF\vert_U \subset T_V $. 
Then the pullback foliation $f^*\sF$ is defined to be the saturation of $\sG$. 
By \cite[Lemma 1.8]{hacon2021birational}, this definition is well-defined. 

Also if $\sG$ is a foliation on $Y$, then we can define the \emph{pushforward foliation} $f_*\sG$ by taking the saturation of the image of the composition 
\[f_*T_\sG\rightarrow f_*T_Y\rightarrow (f_*T_Y)^{**} = T_X.\]
\end{defn}

\begin{defn}[Invariant curves]
Let $(X,\sF)$ be a foliated surface and $U$ be the non-singular locus of $X$. 
A curve $C$ on $X$ is called \emph{$\sF$-invariant} if the inclusion map
\[T_{\sF\vert_U}\vert_C \rw T_U\vert_C\]
factors through $T_C\vert_U$. 
\end{defn}

\begin{defn}[{\cite[Definition I.1.5]{mcquillan2008canonical}}]
Let $(X,\sF,\Delta)$ be a foliated triple and $f : Y\rw X$ be a proper birational morphism. 
For any divisor $E$ on $Y$, we define the \emph{discrepancy} of $(\sF,\Delta)$ along $E$ to be $a(E,\sF,\Delta) = \textnormal{ord}_E(K_{f^*\sF}-f^*(K_\sF+\Delta))$. 
We say  
\[(X,\sF,\Delta) \mbox{ is } \left\{\begin{array}{ll}
\textnormal{terminal} & \mbox{ if } a(E,\sF,\Delta)> 0 \\
\textnormal{canonical} & \mbox{ if } a(E,\sF,\Delta)\geq 0 \\
\textnormal{log terminal} & \mbox{ if } a(E,\sF,\Delta)> -\ve(E) \\
\textnormal{log canonical} & \mbox{ if } a(E,\sF,\Delta)\geq -\ve(E) 
\end{array}\right. 
\mbox{ for every divisor $E$ over $X$} \]
where $\varepsilon(E)$ is defined to be $0$ if $E$ is $f^*\sF$-invariant, and $1$ otherwise. 
\end{defn}

\subsection{Indices on foliated surfaces}
Most definitions in this subsection follow from \cite{brunella2015birational} with some generalizations. 

Let $p\in\textnormal{Sing}(\sF)\backslash\textnormal{Sing}(X)$. That is, $p$ is a smooth point on $X$ but a singular point of the foliation. 
Let $v$ be the vector field around $p$ generating $\sF$. 
Since $p\in\textnormal{Sing}(\sF)$, we have $v(p)=0$. 
Then we can consider the eigenvalues $\lambda_1$, $\lambda_2$ of $(Dv)\vert_p$, which do not depend on the choice of $v$. 

\begin{defn}\label{semi_reduced}
If one of the eigenvalues is non-zero, say $\lambda_2$, then we say $p$ is \emph{semi-reduced} and define the eigenvalue of the foliation $\sF$ at $p$ to be 
\[\lambda := \frac{\lambda_1}{\lambda_2}.\] 
For $\lambda\neq 0$, this definition is unique up to reciprocal $\lambda \sim \frac{1}{\lambda}$. 

If $\lambda = 0$, then $p$ is called a \emph{saddle-node}; otherwise, we say $p$ is \emph{non-degenerate}. 
If $\lambda\not\in\bQ^+$, then $p$ is called a \emph{reduced} singularity of $\sF$. 
\end{defn}

Reduced singularities arise naturally. 
Indeed, blowing up a smooth foliation point will introduce a reduced singularity with $\lambda=-1$. 

\begin{thm}[Seidenberg's theorem]\label{seidenberg}
Given any foliated surface $(X,\sF)$ with $X$ smooth. 
There is a sequence of blowups $\pi : (Y,\sG) \rightarrow (X,\sF)$ such that $(Y,\sG)$ has only reduced singularities.
\end{thm}
\begin{pf}
See \cite{seidenberg1968reduction}, \cite[Appendix]{mattei1980holonomie}, or \cite[Theorem 1.1]{brunella2015birational}. \qed
\end{pf}

\begin{defn}[Separatrices]
Let $p$ be a singular point of $\sF$. 
A \emph{separatrix} of $\sF$ at $p$ is a holomorphic (possibly singular) irreducible $\sF$-invariant curve $C$ on a neighborhood of $p$ which passes through $p$. 
\end{defn}

\begin{thm}\label{separatrix}
Let $\sF$ be a foliation on a normal surface $X$ and $C\subset X$ be a connected, compact $\sF$-invariant curve such that 
\begin{enumerate}
\item all the singularities of $\sF$ on $C$ are reduced. 
\item the intersection matrix of $C$ is negative definite and the dual graph is a tree.
\end{enumerate}
Then there exists at least one point $p\in C\cap\textnormal{Sing}\sF$ and a separatrix through $p$ not contained in $C$. 
\end{thm}
\begin{pf}
See \cite{camacho1988quadratic}, \cite{sebastiani1997existence}, or \cite[Theorem 3.4]{brunella2015birational}. \qed
\end{pf}

\subsubsection{Non-invariant curves}
We first consider the non-invariant curves and define the tangency order for them. 
\begin{defn}
Let $(X,\sF)$ be a foliated surface and $C$ be a \emph{non-invariant} reduced curve. 
Let $p\in C\backslash\textnormal{Sing}(X)$ and $v$ be the vector field generating $\sF$ around $p$. 
Let $f$ be the local defining function of $C$ at $p$.  
We define the \emph{tangency order} of $\sF$ along $C$ at $p$ to be 
\[\textnormal{tang}(\sF,C,p) := \dim_\bC\frac{\cO_{X,p}}{\langle f,v(f)\rangle}.\]
Note that $\mbox{tang}(\sF,C,p)\geq 0$ and is independent of the choices of $v$ and $f$. 
Moreover, if $\sF$ is transverse to $C$ at $p$, then $\textnormal{tang}(\sF,C,p) = 0$. 
Therefore, if $C$ is compact, then we can define 
\[\textnormal{tang}(\sF,C) := \sum_{p\in C}\textnormal{tang}(\sF,C,p).\]
\end{defn}

\begin{prop}[{\cite{brunella1997feuilletages},\cite[Proposition 2.2]{brunella2015birational}}]\label{non_inv}
Let $\sF$ be a foliation on a smooth projective surface $X$, and $C$ be a non-invariant curve on $X$. 
Then we have 
\[K_\sF\cdot C = \textnormal{tang}(\sF,C) - C^2.\]
\end{prop}

\begin{cor}
Let $C$ be a non-invariant curve on a foliated surface $(X,\sF)$. 
If $C$ is contained in the smooth locus of $X$, then we have $(K_\sF+C)\cdot C\geq 0$. 
\end{cor}

\subsubsection{Invariant curves}
Now we study the invariant curves. 
\begin{defn}
Let $(X,\sF)$ be a foliated surface and $C$ be an invariant curve. 
Let $p\in C\backslash\textnormal{Sing}(X)$ and $\omega$ be a $1$-form generating $\sF$ around $p$. 
If $C$ is an invariant curve and $f$ is the local defining function of $C$ at $p$, then we can write \[g\omega = hdf+f\eta\] where $g$ and $h$ are holomorphic functions, $\eta$ is a holomorphic $1$-form, and $h, f$ are relatively prime functions. 

We define the index $\textnormal{Z}(\sF,C,p)$ to be the vanishing order of $\frac{h}{g}\vert_C$ at $p$. 
Also we define the Camacho-Sad index $\textnormal{CS}(\sF,C,p)$ to be the residue of $\frac{-1}{h}\eta\vert_C$ at $p$.  
These two definitions are independent of the choices of $f, g, h, \omega, \eta$. 
(For a reference, see \cite[page 15 in Chapter 2 and page 27 in Chapter 3]{brunella2015birational}.)
\end{defn}

Note that if $p\not\in\textnormal{Sing}(\sF)$, then $\textnormal{Z}(\sF,C,p)=0=\textnormal{CS}(\sF,C,p)$. 
Therefore, if $C$ is compact, then we can define 
\begin{eqnarray*}
\textnormal{Z}(\sF,C) &:=& \sum_{p\in C} \textnormal{Z}(\sF,C,p), \mbox{ and} \\
\textnormal{CS}(\sF,C) &:=& \sum_{p\in C} \textnormal{CS}(\sF,C,p) 
\end{eqnarray*}
where the sums are taken over only finitely many points. 

\begin{defn}
We define the \emph{virtual Euler characteristic} $\chi(C)$ of $C$ on a normal surface $X$ to be $\chi(C) = -K_X\cdot C-C^2$. 
\end{defn}

\begin{thm}[{\cite{brunella1997feuilletages},\cite[Proposition 2.3 and Theorem 3.2]{brunella2015birational}}]\label{CS_formula}
Let $\sF$ be a foliation on a smooth projective surface $X$, and $C$ be an invariant curve on $X$. 
Then we have 
\begin{enumerate}
\item $K_\sF\cdot C = \textnormal{Z}(\sF,C) - \chi(C)$, and 
\item $C^2 = \textnormal{CS}(\sF,C)$. This identity is called Camacho-Sad formula. 
\end{enumerate}
\end{thm}

\begin{lem}\label{lem_ZCS}
Given a foliated surface $(X,\sF)$ and $p$ a reduced singularity. 
\begin{enumerate}
\item If $p$ is non-degenerate, assume $\omega = \lambda y(1+o(1))dx-x(1+o(1))dy$ generates $\sF$ around $p$. 
Then 
\[\textnormal{CS}(\sF,x=0,p) = \frac{1}{\lambda}, \textnormal{CS}(\sF,y=0,p) = \lambda, \mbox{ and } \textnormal{Z}(\sF,x=0,p)=\textnormal{Z}(\sF,y=0,p)=1.\]
\item If $p$ is a saddle-node, assume $\omega = y^{k+1}dx-(x(1+\nu y^k)+yo(k))dy$ generates $\sF$ around $p$ where $k\in\bN$ and $\nu\in\bC$. Then 
\[\textnormal{CS}(\sF,y=0,p) = 0 \mbox{ and } \textnormal{Z}(\sF,y=0,p) = 1.\]
Suppose there exists a weak separatrix, then 
\[\textnormal{CS}(\sF,x=0,p) = \nu \mbox{ and } \textnormal{Z}(\sF,x=0,p) = k+1.\]
\end{enumerate}
\end{lem}
\begin{pf}
This is done by direct computation. 
For a reference, see \cite[page 30-31 in Chapter 3]{brunella2015birational}. 
\qed
\end{pf}

\subsection{Minimal resolutions of foliated surfaces}\label{min_rsln_fol}
\begin{defn}
A morphism $\pi : (Y,\sG) \rightarrow (X,\sF)$ of foliated surfaces is called a \emph{resolution} if $(Y,\sG)$ has only reduced foliation singularities. 

A resolution $\pi : (Y,\sG) \rightarrow (X,\sF)$ of a foliated surface $(X,\sF)$ is called \emph{minimal} if any resolution $\phi : (Z,\sH)\rightarrow (X,\sF)$ of the foliated surface $(X,\sF)$ factors through $\pi$. 
That is, there is a morphism $\psi : (Z,\sH) \rightarrow (Y,\sG)$ with $\sH = \psi^*\sG$ such that $\phi = \pi\circ\psi$.
\end{defn}

\begin{prop}
Any foliated surface $(X,\sF)$ has a unique minimal resolution up to isomorphism. 
\end{prop}
\begin{pf}
Taking a resolution of $X$, and then by Seidenberg's theorem~\ref{seidenberg}, we have a resolution $\pi : (Y,\sG) \rightarrow (X,\sF)$ of the foliated surface $(X,\sF)$. 
After blowing down $\sG$-exceptional curves, we may assume that $(Y,\sG)$ has no $\sG$-exceptional curves. 
We will show that this $\pi$ gives a minimal resolution of the foliated surface $(X,\sF)$. 

Given any resolution $\phi : (Z,\sH) \rightarrow (X,\sF)$ of the foliated surface $(X,\sF)$. 
Let $W$ be the minimal resolution of singularities of $Z\times_XY$, and $\cE$ be the pullback foliation on $W$. 
So we have the following diagram of foliated surfaces 
\[\xymatrix{(W,\cE) \ar[r]^-\alpha \ar[d]_-\beta & (Z,\sH) \ar[d]^-\phi \\ (Y,\sG) \ar[r]_-\pi & (X,\sF).}\]
We may assume $(W,\cE)$ is minimal in the sense that there is no birational morphism of foliated surfaces $\theta : (W,\cE) \rightarrow (W',\cE')$, which is not an isomorphism, 
such that both $\alpha$ and $\beta$ factor through $\theta$. 

Suppose $\alpha$ is not an isomorphism, then there is a (last) $\alpha$-exceptional curve $\widetilde{C}$ on $W$ with $\widetilde{C}^2=-1$, which is also $\cE$-exceptional. 
By the minimality of $(W,\cE)$, $\widetilde{C}$ is not contracted by $\beta$. 
Let $C\subset Y$ be the curve $\beta(\widetilde{C})$. 
Since $p_a(\widetilde{C})=0$, we have that $p_a(C)=0$. 
Note that $C$ is contracted by $\pi$ because $\widetilde{C}$ is contracted by $\phi\circ\alpha = \pi\circ\beta$. 
Hence $C^2\leq -1$. 
Also, $C^2 \geq \widetilde{C}^2 = -1$ and thus, $C^2=-1$. 
Moreover, $\beta$ is isomorphic around $\widetilde{C}$. 
Therefore, $C$ is $\sG$-exceptional. 
But this is impossible since $(Y,\sG)$ has no $\sG$-exceptional curves. 
\qed
\end{pf}

\begin{rmk}
In general, for any minimal resolution $\pi: (Y,\sG) \rightarrow (X,\sF)$ of the foliated surface $(X,\sF)$, the morphism $\pi: Y\rightarrow X$ of surfaces is not the minimal resolution of $X$. 
\end{rmk}

\subsection{Dual graphs}
\begin{defn}
Let $C = \bigcup C_i$ be a collection of proper curves on the smooth locus of a surface $X$. 
Then the (weighted) \emph{dual graph} $\Gamma = \Gamma(C)$ of $C$ is defined as follows:
\begin{enumerate}
\item The vertices of $\Gamma$ are the curves $C_i$. We will use $[C_i]$ to indicate the vertex of $\Gamma$ corresponding to the curve $C_i$.  
\item Two vertices $[C_i]$ and $[C_j]$ for $i\neq j$ are connected with $C_i\cdot C_j$ edges. 
\item The weight $\textnormal{w}([C_i])$ of vertex $[C_i]$ is given by $-C_i^2$.
\end{enumerate}
\end{defn}

\begin{defn}\label{defn}
For any (dual) graph $\Gamma$, we have the following definitions:
\begin{enumerate}
\item A \emph{cycle} is a graph whose vertices and edges can be ordered as $[C_1],\ldots, [C_m]$ and $e_1,\ldots, e_m$ where $m\geq 2$ such that edge $e_i$ connects vertices $[C_i]$ and $[C_{i+1}]$ for $i=1,\ldots, m$ where $[C_{m+1}] := [C_1]$. 
\item We say $\G$ \emph{has simple edges} if any two vertices are connected by at most one edge. 
\item A \emph{tree} is a connected graph which has no cycle as its subgraph. 
\item The \emph{degree}, denoted by $\textnormal{deg}\,[C_i]$, of the vertex $[C_i]$ is the number of edges connecting to $[C_i]$. 
\item The vertex $[C_i]$ is called a \emph{leaf} of $\Gamma$ if the degree of $[C_i]$ is 1. 
\item The vertex $[C_i]$ is called a \emph{fork} of $\Gamma$ if the degree of $[C_i]$ is at least 3. 
\item A tree without forks is called a \emph{chain}. 
\item The connected components of a tree minus a fork are called the \emph{branches} of the fork.  
\item A tree is \emph{star-shaped} if there is exactly one fork. 
\item We define $\D(\G)$ to be the absolute value of the determinant of the intersection matrix for $\G$. 
We set $\D(\emptyset)=1$. 
\end{enumerate}
\end{defn}

By direct computation, we have the following lemma.
\begin{lem}\label{det}
Suppose $\G$ is a tree with simple edges. 
\begin{enumerate}
\item For any vertex $[C]$ of $\G$ of weight $w([C])$, we have 
\[\D(\G) = \textnormal{w}([C])\D\big(\G\bs \{[C]\}\big)-\sum_{i=1}^s\D\big(\G\bs \{[C], [C_i]\}\big)\]
where $C_1,\ldots,C_s$ are all vertices adjacent to $C$. 
\item For any two vertices $[C_i]$ and $[C_j]$ of $\G$, the $(i,j)$-cofactor of the intersection matrix $A$ of $\G$ is 
\[A_{i,j} = (-1)^{n+1}\D\big(\G\bs\{\textnormal{path from } [C_i] \textnormal{ to } [C_j]\}\big)\]
where $n = |\G|$. 
\end{enumerate} 
\end{lem}

\section{Classification}
In this section, we prove Theorem~\ref{main} which gives a full list of all possibilities of the dual graphs of the exceptional divisors on the minimal resolutions of log canonical foliation singularities on surfaces. 

In order to describe the exceptional divisor over a foliation singularity,  we recall the following definitions. 
(See also \cite[Definition 5.1 and 8.1]{brunella2015birational} and \cite[Definition III.0.2 and III.2.3]{mcquillan2008canonical}.)

\begin{defn}\label{-2curve}
Given a foliated surface $(X,\sF)$ with $X$ smooth. 
\begin{enumerate}
\item $E = \bigcup_{i=1}^sE_i$ is called a string if 
\begin{enumerate}
\item each $E_i$ is a smooth rational curve and 
\item $E_i\cdot E_j = 1$ if $|i-j|=1$ and $0$ if $|i-j|\geq 2$. 
\end{enumerate}
\item If, moreover, $E_i^2\leq -2$ for all $i$, then we call $E$ a Hirzebruch-Jung string. 
\item A curve $C\subset X$ is called $\sF$-exceptional if 
\begin{enumerate}
\item $C$ is a $(-1)$-curve. 
\item The contraction of $C$ gives a foliated surface $(X',\sF')$ with only \emph{reduced} singularities. 
\end{enumerate}

\item $C$ is called a $(-1)$-$\sF$-curve (resp. $(-2)$-$\sF$-curve) if 
\begin{enumerate}
\item $C$ is a smooth rational curve. 
\item $\textnormal{Z}(\sF,C)=1$ (resp. $\textnormal{Z}(\sF,C)=2$). 
\end{enumerate}

\item We say $C=\bigcup_{i=1}^sC_i$ is an $\sF$-chain if 
\begin{enumerate}
\item $C$ is a Hirzebruch-Jung string.
\item Each $C_i$ is $\sF$-invariant. 
\item $\textnormal{Sing}(\sF)\cap C$ are all reduced and non-degenerate. 
\item $\textnormal{Z}(\sF, C_1) = 1$ and $\textnormal{Z}(\sF,C_i) = 2$ for all $i\geq 2$. 
\end{enumerate}

\item If an $\sF$-invariant curve $E$ is not contained in a $\sF$-chain $C$ but meets the chain, 
then we call $E$ the tail of the chain $C$. 

\item $C$ is called a bad tail if 
\begin{enumerate}
\item $C$ is a smooth rational curve with $\mbox{Z}(\sF,C) = 3$ and $C^2\leq -2$. 
\item $C$ intersects two $(-1)$-$\sF$-curves whose self-intersections are both $-2$. 
\end{enumerate}
\end{enumerate}
\end{defn}

Before proving the main theorem, we recall the following well-known lemma. 
\begin{lem}\label{key_lem}
Let $C = \bigcup_iC_i$ be a set of proper curves on a smooth surface. 
Assume that the intersection matrix $(C_i\cdot C_j)_{i,j}$ is negative definite. 
Let $A = \sum a_iC_i$ be an $\bR$-linear combination of the curves $C_i$'s. 
If $A\cdot C_j\geq 0$ for all $j$, then 
\begin{enumerate}
\item $a_i\leq 0$ for all $i$. 
\item If $C$ is connected, then either $a_i=0$ for all $i$ or $a_i<0$ for all $i$. 
Moreover, if $A\cdot C_j> 0$ for some $j$, then $a_i <0$ for all $i$. 
\end{enumerate}
\end{lem}
\begin{pf}
For a reference, see~\cite[Lemma 3.41]{kollar1998birational}. 
\end{pf}

\begin{rmk}
The exceptional divisor $E$ of a resolution of a singularity is connected by Zariski's main theorem. 
Moreover, it is well-known that the intersection matrix of $E$ is negative definite. (For a reference, see \cite[Lemma 3.40]{kollar1998birational}.)
\end{rmk}

\begin{thm}\label{main}
Let $(X,\sF,p)$ be a germ of a foliated surface. 
Assume that $p$ is a log canonical singularity of $\sF$. 
Let $\pi : (Y,\sG) \rightarrow (X,\sF)$ be the minimal resolution (see subsection~\ref{min_rsln_fol}) for $(X,\sF,p)$ with exceptional divisors $E = \cup E_i$. 
Then $E$ belongs to one of the following types: 
\begin{enumerate}
\item A $\sG$-chain. 
\item A chain of three invariant curves $E_1\cup E_2\cup E_3$ where $E_1$ and $E_3$ are $(-1)$-$\sG$-curves with self-intersection $-2$ and $E_2$ is a bad tail.
\item A chain of $(-2)$-$\sG$-curves. 
\item A dihedral singularity. 
More precisely, two $(-1)$-$\sG$-curves with self-intersection $-2$ joined by a bad tail which itself connects to a chain of $(-2)$-$\sG$-curves. 
\item An elliptic Gorenstein leaf. 
\item A chain $E = \bigcup_{i=1}^rE_i$ with exactly one non-invariant curve $E_\ell$ with $1\leq \ell\leq r$. Moreover, $E_\ell$ has tangency order zero and $\bigcup_{i=1}^{\ell-1}E_i$ and $\bigcup_{i=\ell+1}^{r}E_i$ are $\sG$-chains. 
\item The dual graph is star-shaped with a non-invariant center $[E_0]$. Moreover, $E_0$ has tangency order zero, all branches are $\sG$-chains, and all first curves of $\sG$-chains have intersection number one with $E_0$. 
\end{enumerate}
Note that type (1) is terminal, and types (1) - (5) are canonical. 
\end{thm}

\begin{pf}
Let $\Gamma$ be the dual graph of exceptional divisors $E=\cup E_i$, which is connected by Zariski's main theorem. 
Let $\Delta$ be the sum over all non-invariant exceptional curves $E_i$. 
Then we write $K_\sG+\Delta = \pi^*K_\sF + \sum a_iE_i$ where $a_i\geq 0$ since $p$ is a log canonical singularity of $\sF$. 

We divide the proof into several steps. 
\begin{enumerate}
\item\label{inv_deg} \emph{Claim.} For any vertex $[E_i]$ with $E_i$ invariant, we have $(K_\sG+\Delta)\cdot E_i \geq \deg\,[E_i] -2$. 
\begin{pf}
Let $d = \deg\,[E_i]$. 
Suppose there are exactly $m$ edges of $[E_i]$ connecting to the vertices corresponding to \emph{invariant} curves. 
Since the intersections of two invariant curves are reduced singularities and at most two separatrices pass through any reduced singularity, there are at least $m$ distinct foliation singularities on $E_i$. 
Thus, we have that $\textnormal{Z}(\sG,E_i)\geq m$. 

The other $d-m$ edges of $[E_i]$ connect to the vertices corresponding to \emph{non-invariant} curves, which are in the support of $\Delta$. 
Thus, we have that $\Delta\cdot E_i = d-m$ by the definition of edges. 
Therefore, by Theorem~\ref{CS_formula}, we have 
\[(K_\sG+\Delta)\cdot E_i = \textnormal{Z}(\sG,E_i)+2p_a(E_i)-2+\Delta\cdot E_i\geq m -2 + (d-m) = d-2.\]
\end{pf}

\item\label{one_vertex} 
\emph{Claim.} Suppose $\Gamma = \{[E_1]\}$ is a one-vertex graph. Then we have the following two possibilities. 
\begin{enumerate}
\item If $E_1$ is non-invariant, then the tangency order of $E_1$ is zero. 
\item If $E_1$ is invariant, then $E_1$ is a $\sG$-chain, a $(-2)$-$\sG$-curve, or a rational curve with only one node. 
\end{enumerate}

\begin{pf}
If $E_1$ is non-invariant, then we have $\Delta=E_1$ and, by Theorem~\ref{non_inv}, 
\[(K_\sG+\Delta)\cdot E_1 = (K_\sG+E_1)\cdot E_1 = \textnormal{tang}(\sG,E_1)\geq 0. \]
By Lemma~\ref{key_lem}, we have $a_1\leq 0$ and $a_1=0$ if and only if $(K_\sG+\Delta)\cdot E_1=0$. 
Since $a_1\geq 0$, we have $a_1=0$ and hence 
\[\textnormal{tang}(\sG,E_1) = (K_\sG+E_1)\cdot E_1 = (K_\sG+\Delta)\cdot E_1 = \pi^*K_\sF\cdot E_1=0.\]

If $E_1$ is invariant, then we hvae $\Delta = 0$ and, by Theorem~\ref{CS_formula}, 
\[K_\sG\cdot E_1 = \textnormal{Z}(\sG,E_1)+2p_a(E_1)-2 \geq -2.\]
\begin{itemize}
\item If $K_\sG\cdot E_1=-2$, then $\textnormal{Z}(\sG,E_1)=0$ and $p_a(E_1)=0$. 
So $E_1$ is smooth with no foliation singularity, which is impossible by Theorem~\ref{separatrix}.  
\item If $K_\sG\cdot E_1=-1$, then $\mbox{Z}(\sG,E_1)=1$ and $p_a(E_1)=0$. 
So $E_1$ is a smooth rational curve with exactly one reduced singularity. 
Since $\pi$ is minimal, we have $E_1^2\leq -2$, and thus $E_1$ is a $\sG$-chain. 
\item If $K_\sG\cdot E_1\geq 0$, then by Lemma~\ref{key_lem}, we have $a_1\leq 0$ and $a_1=0$ if and only if $K_\sG\cdot E_1=0$. 
Since $a_1\geq 0$, we have $a_1=0$ and $K_\sG\cdot E_1=0$. 
Thus we have $p_a(E_1) = 0$ or $1$. 
If $p_a(E_1)=0$, then $\textnormal{Z}(\sG,E_1) = 2$, which is a $(-2)$-$\sG$-curve. 
If $p_a(E_1)=1$, then $\textnormal{Z}(\sG,E_1) = 0$. 
So $E_1$ is not smooth, otherwise $E_1^2 = 0$ by Theorem~\ref{CS_formula}, which is impossible. 
Thus, $E_1$ is a rational curve with only one node. 
\end{itemize}
\end{pf}

\item\label{lower_bound} 
From now on, we assume that $\Gamma$ has at least two vertices. 
\begin{claim}
For any vertices $[E_i]$, we have 
$(K_\sG+\Delta)\cdot E_i\geq -1$ and the equality holds if and only if $E_i$ is invariant and $[E_i]$ is a leaf. 
\end{claim}
\begin{pf}
Since $\Gamma$ is connected with at least two vertices, we have $\deg\,[E_i]\geq 1$ for any vertex $[E_i]$. If $E_i$ is invariant, then we have the required inequality by step~(\ref{inv_deg}). 

If $E_i$ is non-invariant, then we have $E_i\subset\Delta$ and \[(K_\sG+\Delta)\cdot E_i\geq (K_\sG+E_i)\cdot E_i = \textnormal{tang}(\sG,E_i)\geq 0 \]
by Theorem~\ref{non_inv}. 
\end{pf}

\item We call $[E_i]$ a \emph{special leaf} if $E_i$ is invariant with $(K_\sG+\Delta)\cdot E_i=-1$. 
We will denote $L := E_i$ if $[E_i]$ is a special leaf. 

Note that, by Theorem~\ref{CS_formula}, we have 
\[-1 = (K_\sG+\Delta)\cdot L = \textnormal{Z}(\sG,L) + 2p_a(L)-2 + \Delta\cdot L\geq -2.\]
Hence, we have $p_a(L)=0$, and thus $L$ is a smooth rational curve. 
Since $L^2\leq -1$, we have  $\textnormal{Z}(\sG,C) \geq 1$, and hence $\textnormal{Z}(\sG,C) = 1$ and $\Delta\cdot L =0$. 
Then by the minimality of $\pi$, we have $L^2\leq -2$. 
Thus $L$ itself is a $\sG$-chain. 
Therefore, we consider the chain $C_L$ of maximal length $N_L$ starting from $L$, which is a $\sG$-chain and disjoint from $\Delta$. 

\item\label{tail} 
\emph{Claim.} Fix a special leaf $[L]$. Other than those curves in $C_L$, there is at most one curve in $E$ 
connecting to the last curve $C_{L,N_L}$ of $C_L = \bigcup_{i=1}^{N_L}C_{L,i}$ with $C_{L,1}=L$ and $C_{L,j}\cdot C_{L,j+1}=1$ for $j=1,\ldots, N_L-1$. 
Precisely, if the dual graph of $C_L$ is not $\Gamma$, then such curve exists and is called the tail of $C_L$, denoted by $C_{L,N_L+1}$.  
Moreover, it is invariant with 
\[\textnormal{Z}(\sG,C_{L,N_L+1}) \geq 2 \mbox{ and } (K_\sG+\Delta)\cdot C_{L,N_L+1}\geq 1.\] 

\begin{pf}
If the dual graph of $C_L$ is $\Gamma$, then the statement is clear. 
If not, then such vertex $[V]$ exists because $\Gamma$ is connected. 
Since $C_L$ is disjoint from $\Delta$, so the curve $V$ must be invariant. 
Since $C_L$ is a $\sG$-chain, there is precisely one such curve $V$. 
Besides, this curve $V$ must connect to the last curve $C_{L,N_L}$. 
So we denote it by $C_{L,N_L+1}$. 
Note that $\textnormal{Z}(\sG,C_{L,N_L+1}) \geq 1$ since $C_{L,N_L+1}$ intersects $C_L$. 
\begin{enumerate}
\item If $\textnormal{Z}(\sG,C_{L,N_L+1}) = 1$, then $\bigcup_{j=1}^{N_L+1}C_{L,j}$ contradicts to Theorem~\ref{separatrix}. 
\item If $\textnormal{Z}(\sG,C_{L,N_L+1}) = 2$, then $\Delta\cdot C_{L,N_L+1} \geq 1$ by the maximality of $C_L$. 
Thus, by Theorem~\ref{CS_formula}, we have 
\[(K_\sG+\Delta)\cdot C_{L,N_L+1} = \textnormal{Z}(\sG,C_{L,N_L+1}) - 2 + \Delta\cdot C_{L,N_L+1}\geq 1.\] 
\item If $\textnormal{Z}(\sG,C_{L,N_L+1}) \geq 3$, then $(K_\sG+\Delta)\cdot C_{L,N_L+1}\geq 1$. 
\end{enumerate} 
\end{pf}

\item\label{not_chain} From now on, we assume that $\Gamma$ is not a dual graph of a $\sG$-chain. 
So by step~(\ref{tail}), $C_{L,N_L+1}$ always exists for any special leaf $[L]$. 

\item\label{disjoint} \emph{Claim.} For two distinct special leaves $[L]$ and $[L']$, two chains $C_L$ and $C_{L'}$ are disjoint. 
\begin{pf}
If not, then two $\sG$-chains intersect at the last curves, and $C_L\cup C_{L'}$ forms a chain which contradicts Theorem~\ref{separatrix}. 
\end{pf}

\item Now for any special leaf $[L]$, we define a divisor 
\[D_L = \sum_{j=1}^{N_L} \frac{N_L+1-j}{N_L+1}C_{L,j}.\]
Let $D=\sum D_L$ where the sum is over all special leaves $[L]$. 

\item\label{zero} 
\emph{Claim.} There is a vertex $[V]$ of $\Gamma$ such that the coefficient of $V$ in $D$ is zero. 
\begin{pf}
Suppose not, then the support of $D$, which is a union of $C_L$'s where $[L]$'s are special leaves, is the same as the support of $E$. 
Since $\Gamma$ is connected and any two distinct $C_L$ and $C_{L'}$ are disjoint by step~(\ref{disjoint}), we conclude that $\Gamma$ is the dual graph of $C_L$ for one special leaf $[L]$. 
Therefore, $\Gamma$ is the dual graph of a $\sG$-chain, which contradicts step~(\ref{not_chain}). 
\end{pf}

\item 
\emph{Claim.} We have $D\cdot L \leq -1 = (K_\sG+\Delta)\cdot L$ for any special leaf $[L]$. 
Moreover, $D\cdot L=-1$ if and only if $L^2=-2$. 

\begin{pf}
If $N_L=1$, then 
\[D\cdot L = \frac{1}{2}L^2\leq -1\]
and the equality holds if and only if $L^2=-2$. 

If $N_L\geq 2$, then 
\[D\cdot L = \frac{N_L}{N_L+1}L^2+\frac{N_L-1}{N_L+1}C_{L,2}\cdot L \leq \frac{-2N_L}{N_L+1}+\frac{N_L-1}{N_L+1} = -1\] 
and the equality holds if and only if $L^2 = -2$. 
\end{pf}

\item \emph{Claim.} For any special leaf $[L]$, we have 
\[D\cdot C_{L,j} \leq 0 = (K_\sG+\Delta)\cdot C_{L,j}\] 
for $2\leq j\leq N_L$. 

\begin{pf}
Since $C_{L,j}$ is a $(-2)$-$\sG$-curve and $C_L$ is disjoint from $\Delta$, we have $K_\sG\cdot C_{L,j} = 0$ and $\Delta\cdot C_{L,j} = 0$. 
Therefore, we have $(K_\sG+\Delta)\cdot C_{L,j}=0$. 

For $D\cdot C_{L,j}$, we have 
\begin{eqnarray*}
D\cdot C_{L,j} &=& \frac{N_L+2-j}{N_L+1}+\frac{N_L+1-j}{N_L+1}C_{L,j}^2+\frac{N_L-j}{N_L+1} \\
&\leq & \frac{N_L+2-j}{N_L+1}+\frac{N_L+1-j}{N_L+1}(-2)+\frac{N_L-j}{N_L+1} = 0.
\end{eqnarray*} 
\end{pf}

\item \emph{Claim.} For any fork $[F]$, we have 
\[D\cdot F = \sum \frac{1}{N_L+1}\] 
where the sum is over all special leaves $[L]$ such that $C_{L,N_L+1}=F$. 
Also we have $D\cdot F\leq (K_\sG+\Delta)\cdot F$.

\begin{pf}
It is clear that 
\[D\cdot F = \sum\frac{1}{N_L+1}C_{L,N_L}\cdot F = \sum\frac{1}{N_L+1}\] 
where the sum is over all special leaves $[L]$ such that $C_{L,N_L+1}=F$. 
\begin{itemize}
\item If $F$ is non-invariant, then $D\cdot F = 0$ since $C_{L,N_L+1}$ is always invariant by step~(\ref{tail}). 
Also, by Theorem~\ref{non_inv}, we have $(K_\sG+\Delta)\cdot F\geq (K_\sG+F)\cdot F\geq 0$. 
\item If $F$ is invariant, then $D\cdot F \leq m/2$ where $m$ is the number of special leaves $[L]$ such that $C_{L,N_L+1} = F$. 
Also, by Theorem~\ref{CS_formula}, we have 
\[(K_\sG+\Delta)\cdot F\geq K_\sG\cdot F \geq \textnormal{Z}(\sG,F)-2 \geq m-2.\] 
Moreover, by step~(\ref{inv_deg}), we have $(K_\sG+\Delta)\cdot F\geq \deg\,[F]-2\geq 1$. 
Hence, we have 
\[(K_\sG+\Delta)\cdot F\geq\max\{m-2,1\} \geq \frac{m}{2} \geq D\cdot F.\]
\end{itemize}
\end{pf}

\item\label{requirement} 
\emph{Claim.} For any vertex $[V]$ of $\Gamma$, we have $D\cdot V\leq (K_\sG+\Delta)\cdot V$. 

\begin{pf}
The only vertices we have not checked yet are vertices $[V]$ of degree $\leq 2$, which are not in any $C_L$ where $[L]$ is a special leaf. 
If $V = C_{L,N_L+1}$ for some special leaf $[L]$, then we have 
\[D\cdot V \leq \frac{1}{N_{L_1}+1} + \frac{1}{N_{L_2}+1} \leq 1\]
since $\deg\,[V] \leq 2$ and $N_L\geq 1$ for all special leaves $[L]$. 
Also, by step~(\ref{tail}), we have that $(K_\sG+\Delta)\cdot V\geq 1$. 

If $V$ is not $C_{L,N_L+1}$ for any special leaf $[L]$, then $D\cdot V=0$. 
Thus, by step~(\ref{lower_bound}), we have $(K_\sG+\Delta)\cdot V\geq 0$. 
\end{pf}

\item\label{equal} 
\emph{Claim.} $(K_\sG+\Delta)\cdot V = D\cdot V$ for all vertices $[V]$ of $\Gamma$. 

\begin{pf}
Suppose there is a vertex $[V]$ such that $(K_\sG+\Delta)\cdot V \neq D\cdot V$. 
By steps~(\ref{requirement}), (\ref{zero}), and Lemma~\ref{key_lem}, we have that $a_i < 0$ for some $i$, which gives a contradiction. 
\end{pf}

\item\label{bad_tail} 
\emph{Claim.} For any special leaf $[L]$, we have $N_L=1$ and $C_{L,2}$ must be a bad tail which has 
\[\deg\,[C_{L,2}]=2 \mbox{ or } 3\] and is disjoint from $\Delta$. 
\begin{pf}
Let $V := C_{L,N_L+1}$ and $L_1,\ldots, L_s$ be all special leaves such that $C_{L_i,N_{L_i}+1}=V$ for $i=1,\ldots, s$. So we have $\textnormal{Z}(\sG,V)\geq s$. 

If $\textnormal{Z}(\sG,V) = s$, then $V\cup\bigcup_{i=1}^s C_{L_i}$ contradicts to Theorem~\ref{separatrix}. 
Thus, we have 
\begin{equation}\label{s+1}
\textnormal{Z}(\sG,V) \geq s+1
\end{equation} 
and then $(K_\sG+\Delta)\cdot V \geq s-1$. 

For the case when $s=1$, we have $D\cdot V = \frac{1}{N_L+1}\leq\frac{1}{2}$ and $(K_\sG+\Delta)\cdot V \geq 1$ by step~(\ref{tail}), which is a contradiction to step~(\ref{equal}). 

For the case when $s\geq 2$, we have 
\[s-1\leq (K_\sG+\Delta)\cdot V = D\cdot V = \sum_{i=1}^s\frac{1}{n_{L_i}+1} \leq \frac{s}{2}.\] 
This inequality holds only when $s=2$, $N_{L_1}=N_{L_2}=1$, and $(K_\sG+\Delta)\cdot V = 1$. 
Thus, by step~(\ref{inv_deg}), we have $\deg\,[V]\leq 3$ and hence $\deg\,[V]=2$ or $3$. 
Notice that, by Theorem~\ref{CS_formula}, we have  
\[1 = (K_\sG+\Delta)\cdot V \geq \textnormal{Z}(\sG,V) -2 + \Delta\cdot V \geq 3-2 = 1.\]
Therefore, all inequalities are equalities, that is, $\Delta\cdot V=0$, $\textnormal{Z}(K_\sG,V)=3$, and $p_a(V) = 0$.
\end{pf}

\item\label{non_inv_fork} 
\emph{Claim.} Let $[F]$ be a fork. Then $F$ is either a bad tail or a non-invariant curve with tangency order zero. 
\begin{pf}
Assume that $F$ is not a bad tail. 
Because of step~(\ref{bad_tail}), we have $F\neq C_{L,N_L+1}$ for any special leaf $[L]$. 
Therefore, we have $D\cdot F=0$. 

If $F$ is invariant, then $(K_\sG+\Delta)\cdot F\geq 1$ by step~(\ref{inv_deg}), which is a contradiction to step~(\ref{equal}). 

Hence, $F$ is non-invariant with tangency order zero since  
\[0 = D\cdot F = (K_\sG+\Delta)\cdot F \geq (K_\sG+F)\cdot F = \textnormal{tang}(\sG,F)\geq 0.\]
\end{pf}

\item\label{deg2}
\emph{Claim.} All vertices $[V]$ of degree $\leq 2$ belong to one of the following types: 
\begin{enumerate}
\item A non-invariant curve $V$ with tangency order zero. 
\item A special leaf. 
\item A bad tail $V$.
\item An invariant curve $V$ with $p_a(V)=0$, $\textnormal{Z}(\sG,V) = 1$, and $\Delta\cdot V=1$. 
\item An invariant curve $V$ with $p_a(V)=0$, $\textnormal{Z}(\sG,V) = 2$, and $\Delta\cdot V=0$. 
\end{enumerate} 

\begin{pf}
If $V$ is non-invariant, then the tangency order is zero since 
\[0 = D\cdot V = (K_\sG+\Delta)\cdot V \geq (K_\sG+V)\cdot V = \textnormal{tang}(\sG,V)\geq 0.\] 

Now suppose $V$ is invariant but not a bad tail, and $[V]$ is not a special leaf. 
Then, by step~(\ref{bad_tail}), we have $D\cdot V=0$. 
Thus we have $(K_\sG+\Delta)\cdot V = 0$ by step~(\ref{equal}). 
Moreover, $\mbox{Z}(\sG,V)\geq 1$,  otherwise $V^2=0$ by Theorem~\ref{CS_formula}, which is impossible. 
So \[0=(K_\sG+\Delta)\cdot V = \textnormal{Z}(\sG,V)+2p_a(V)-2+\Delta\cdot V\geq 2p_a(V)-1.\]
Hence $p_a(V) = 0$ and $\textnormal{Z}(\sG,V)+\Delta\cdot V = 2$. 
Therefore we have either $\textnormal{Z}(\sG,V) = 1$ and $\Delta\cdot V=1$ or $\textnormal{Z}(\sG,V) = 2$ and $\Delta\cdot V=0$. 
\end{pf}

\item\label{one_fork} 
\emph{Claim.} There is at most one fork. 

\begin{pf}
Suppose there are at least two forks. 
Let $[F_1]$ be one of the forks. 
If $F_1$ is invariant, then it is a bad tail by step~(\ref{non_inv_fork}), and the degree of $[F_1]$ must be three by step~(\ref{bad_tail}). 
Also, two of the three branches of $[F_1]$ are special leaves $[L_1]$ and $[L_2]$. 

Let $[F_2]$ be the fork connecting to $[F_1]$ by a chain of vertices $[V_1],\ldots, [V_r]$ of degree two where $V_1$ intersects $F_1$. Let $V_{r+1}=F_2$. 

If $V_i$ is a non-invariant curve for some $1\leq i\leq r+1$, then we put $i_0$ be the minimum of $i$ such that $V_i$ is non-invariant. 
By step~(\ref{bad_tail}), we have $i_0\geq 2$. 
Also by step~(\ref{deg2}), $[V_j]$ is of type (e) for $1\leq j\leq i_0-2$ and of type (d) for $j = i_0-1$. 
Thus, $L_1\cup L_2\cup F_1 \cup\bigcup_{j=1}^{i_0-1}V_j$ contradicts to Theorem~\ref{separatrix}. 

Therefore, all $C_i$ are invariant for $1\leq i\leq r+1$. 
In particular, $F_2$ is an invariant fork, which is a bad tail by step~(\ref{non_inv_fork}). 
However, $\Gamma$ contradicts Theorem~\ref{separatrix}. 

Hence we may assume that all forks are non-invariant. 
Then there are two (non-invariant) forks connected by a chain of vertices $[V_1],\ldots, [V_r]$ of degree two. 
Let $[V_0]$ and $[V_{r+1}]$ be these two forks. 
Thus by step~(\ref{deg2}), we have that $[V_1]$ and $[V_r]$ are of type (d) 
and $[V_j]$ is of type (e) for $2\leq j\leq r-1$. 
However, $\bigcup_{i=1}^r V_i$ contradicts to Theorem~\ref{separatrix}. 
\end{pf}

\item\label{list_no_fork}
\emph{Claim.} If $\Gamma = \Gamma(E)$ has no fork, then $E$ belongs to one of the following types: 
\begin{enumerate}
\item A $\sG$-chain. 
\item A chain of three invariant curves $E_1\cup E_2\cup E_3$ where $E_1, E_3$ are special leaves of self-intersection $-2$ and $E_2$ is a bad tail.
\item A chain of $(-2)$-$\sG$-curves. 
\item A chain $E = \bigcup_{i=1}^rE_i$ with exactly one non-invariant curve $E_\ell$ with $1\leq \ell\leq r$. 
Moreover, $E_\ell$ has tangency order zero and $\bigcup_{i=1}^{\ell-1}E_i$ and $\bigcup_{i=\ell+1}^{r}E_i$ are $\sG$-chains. 
\item An elliptic Gorenstein leaf. 
\end{enumerate}
\begin{pf}
Note that $\Gamma$ is either a chain or a cycle since $\Gamma$ has no fork. 
Assume that $\Gamma$ is a chain but not the dual graph of a $\sG$-chain. 
If $\Gamma$ has a special leaf, then by step~(\ref{bad_tail}), the tail of the special leaf is actually a bad tail, and the degree of the tail must be two since there is no fork in $\Gamma$. 
So $E$ is of type (b). 

If $\Gamma$ has no special leaf, then $D=0$ and $(K_\sG+\Delta)\cdot V = 0$ for all vertices $[V]$ of $\Gamma$ by step~(\ref{equal}). 
Suppose there are more than two vertices $[V]$ with non-invariant $V$. 
Let $[V_1]$ and $[V_2]$ be two such vertices and be connected by a chain of vertices $[C_1],\ldots, [C_r]$ with $C_i$ invariant for $i=1,\ldots, r$. 
Then by step~(\ref{deg2}), we have that $[C_1]$ and $[C_r]$ are of type (d) and $[C_j]$ is of type (e) for $2\leq j\leq r-1$. 
However, $\bigcup_{i=1}^r C_i$ contradicts to Theorem~\ref{separatrix}. 
Hence $E$ is of type (c) if there is no non-invariant curve and type (d) if there is precisely one non-invariant curve. 

When $\Gamma$ is a cycle, there is no non-invariant curve by the same argument as above, and hence $\Gamma$ is an elliptic Gorenstein leaf. 
\end{pf}

\item \emph{Claim.} If $\Gamma = \Gamma(E)$ has exactly one fork, then $\Gamma$ belongs to one of the following types:
\begin{enumerate}
\item A graph of $D_n$ type. 
More precisely, two $(-1)$-$\sG$-curves with self-intersections $-2$ joined by a bad tail which itself connects to a chain of $(-2)$-$\sG$-curves. 
\item A star-shaped graph with non-invariant center $E_0$. 
Moreover, $E_0$ has tangency order zero, all branches are $\sG$-chains, and all first curves of $\sG$-chains have intersection number one with $E_0$. 
\end{enumerate}
\begin{pf}
Let $[F]$ be the fork. 
If $F$ is invariant, then $F$ must be a bad tail by step~(\ref{non_inv_fork}). 
Let $[V]$ be a vertex other than $[F]$ and its associated two special leaves $[L_1]$ and $[L_2]$. 
Note that $D\cdot V=0$ and therefore we have $(K_\sG+\Delta)\cdot V=0$ by step~(\ref{equal}). 
Note that the graph $\Gamma' := \Gamma\backslash\{[F],[L_1],[L_2]\}$ is a chain by step~(\ref{list_no_fork}). 
Moreover, $\Gamma'$ must be the graph of some curves belonging to one of the types in step~(\ref{list_no_fork}). 
Since $(K_\sG+\Delta).V=0$ for all vertices $[V]$ in $\Gamma'$, the type (a) and (b) in step~(\ref{list_no_fork}) is impossible. Note that $\Gamma'$ is a chain, so the type (e) in step~(\ref{list_no_fork}) is also impossible. 
If $\Gamma'$ is of type (d) in step~(\ref{list_no_fork}), we put $E_1$ be the curve intersecting $F$ and then $\ell\geq 2$ where $E_\ell$ is non-invariant. 
However, $L_1\cup L_2\cup F \cup\bigcup_{i=1}^{\ell-1}E_i$ contradicts to Theorem~\ref{separatrix}. 
Thus, $\Gamma'$ is of type (c) in step~(\ref{list_no_fork}), which gives $\Gamma$ is of type (a). 

Now if $F$ is non-invariant, then there is no special leaf, otherwise, by step~(\ref{bad_tail}), there is a bad tail $B$ with $\deg\,[B]=3$, which is impossible by step~(\ref{one_fork}). 
Thus we have $D=0$, and by step~(\ref{equal}), $(K_\sG+\Delta)\cdot V = D\cdot V=0$ for all vertices $[V]$. 
Moreover, every branch of $[F]$ belongs to one of the types in step~(\ref{list_no_fork}). 

Let $B = \{[B_i]\}_{i=1}^r$ be one of the branches of $[F]$ with $B_1$ intersecting $F$. 
Let $\Delta_B$ be the sum over all non-invariant curves among $B_i$'s. 
Note that 
\[(K_\sG+\Delta_B)\cdot B_1 = (K_\sG+\Delta)\cdot B_1 - F\cdot B_1 \leq -1.\] 
By steps~(\ref{one_vertex}) and (\ref{list_no_fork}), we have that $B_1$ is a $(-1)$-$\sG$-curve and $F\cdot B_1=1$. 
For any $j\geq 2$, we have $(K_\sG+\Delta_B)\cdot B_j = (K_\sG+\Delta)\cdot B_j=0$. 
Since $B$ is a chain and belongs to one of the types in step~(\ref{list_no_fork}), we have that $\bigcup_{i=1}^rB_i$ is a $\sG$-chain. Hence, $\Gamma$ is of type (b). 
\end{pf}\qed
\end{enumerate}
\end{pf}

\section{Foliated minimal log discrepancy}
In this section, we define the foliated minimal log discrepancy and show some of its properties.
\begin{defn}
Let $(X,\sF,\Delta)$ be a foliated triple. 
For any divisor $E$ over $X$, we define the \emph{foliated log discrepancy} of $(X,\sF,\Delta)$ along $E$ to be $a(E,\sF,\Delta) + \ve(E)$ where $\ve(E)=0$ if $E$ is invariant under the pullback foliation and $\ve(E)=1$ otherwise.
\end{defn}

\begin{defn}\label{foliated_discrepancy}
Given $(X,\sF,\Delta)$ a foliated triple. 
Let $(Y,\sG)$ be the minimal resolution of $(X,\sF)$. 
We define the \emph{foliated minimal log discrepancy} of $(X,\sF,\Delta)$ as 
\[\textnormal{mld}(\sF,\Delta) := \inf\{a(E,\sF,\Delta)+\ve(E)\neq 0\vert\,\mbox{$E$ is a divisor over $X$}\}\]
and the \emph{partial log discrepancy} as 
\[\textnormal{pld}(\sF,\Delta) := \min\{a(E,\sF,\Delta)+\ve(E)\neq 0\vert\,\mbox{$E$ is a divisor on $Y$}\}.\] 

Also, for any fixed $x\in X$, we define 
\begin{eqnarray*}
\textnormal{mld}_x(\sF,\Delta) &:=& \inf\{a(E,\sF,\Delta)+\ve(E)\neq 0\vert\,\mbox{the center of $E$ on $X$ is $x$}\} \mbox{ and} \\
\textnormal{pld}_x(\sF,\Delta) &:=& \min\{a(E,\sF,\Delta)+\ve(E)\neq 0\vert\,\mbox{$E$ is a divisor on $Y$ over $x$}.\}
\end{eqnarray*}
From now on, we make the convention that $\min\emptyset = \inf\emptyset = 0$. 
\end{defn}
\begin{rmk}
By Corollary~\ref{mld_min}, if $\textnormal{mld}_x(\sF,\D)\geq 0$, then it is indeed a minimum, that is, there is a divisor $E$ over $X$ that computes the $\textnormal{mld}_x(\sF,\D)$. 
\end{rmk}

\begin{prop}
If $\textnormal{mld}(\sF,\Delta)<0$, then $\textnormal{mld}(\sF,\Delta)=-\infty$.
\end{prop}
\begin{pf}
Suppose $\textnormal{mld}(\sF,\Delta)<0$, then there is a divisor $E$ over $X$ such that \[a(E,\sF,\Delta)+\ve(E)<0.\]

If $E$ is non-invariant, that is $\ve(E)=1$, then blowing up a general point $p$ on $E$ introduces an invariant exceptional divisor $E'$ with $a(E',\sF,\Delta) +\ve(E') \leq a(E,\sF,\Delta)+1<0$. 
Therefore, we may assume that $E$ is invariant. 

Let $a(E,\sF,\Delta) = -c$ for some positive number $c$. 
Now blowing up a general point $p$ on $E$ introduces an invariant exceptional divisor $E_1$ with 
\[a(E_1,\sF,\Delta) \leq 1+a(E,\sF,\Delta) = 1-c.\]
Notice that the proper transform of the support of $\Delta$ doesn't contain the intersection of $E_1$ and the proper transform of $E$. 

Next we blow up the intersection of $E_1$ and the proper transform of $E$, then we have an invariant exceptional divisor $E_2$ with \[a(E_2,\sF,\Delta) = a(E,\sF,\Delta) + a(E_1,\sF,\Delta) \leq 1+ 2a(E,\sF,\Delta) = 1-2c.\]
Then we blow up the intersection of the proper transform of $E$ and $E_2$ to have an invariant exceptional divisor $E_3$ with 
\[a(E_3,\sF,\Delta) \leq 1+ 3a(E,\sF,\Delta) = 1-3c.\]
Repeating the process, by inductively, we have an invariant exceptional divisor $E_n$ over $X$ with $a(E_n,\sF,\D)\leq 1-nc$. 
This shows that $\textnormal{mld}(\sF,\Delta)=-\infty$. 
\qed
\end{pf}

\begin{prop}\label{mld_sm_model}
Let $(X,\sF,\Delta)$ be a foliated triple with only log canonical foliation singularities, that is $\textnormal{mld}(\sF,\Delta)\geq 0$. 
Assume that $\Delta = \sum a_iD_i$ is an $\bR$-divisor which has the simple normal crossing support where $a_i\leq 1$. 
Suppose $X$ is smooth, $\sF$ has only reduced singularities, and any separatrix $C$ through a non-smooth foliation point on non-invariant $D_i$ has $C\cdot D_i=1$, then we have 
\[\textnormal{mld}(\sF,\Delta) = \min\left\{\min_{(i,j)\in S_1}\{1-a_i-a_j\},\min_{(i,j)\in S_0}\{-a_i-a_j\},\min_{i\in I}\{1-a_i\},\min_{i\not\in I}\{-a_i\}, 1\right\}\]
where 
\begin{eqnarray*}
S_1 &=& \{(i,j)\,\vert\, i\neq j, D_i\cap D_j\neq\emptyset \mbox{, and } D_i\cap D_j \mbox{ supports on smooth foliation points}\}, \\
S_0 &=& \{(i,j)\,\vert\, i\neq j \mbox{ and } D_i\cap D_j\neq\emptyset\} \bs S_1\mbox{, and} \\
I &=& \{i\,\vert\, D_i \mbox{ is non-invariant and has only smooth foliation points}\}. 
\end{eqnarray*}
\end{prop}
\begin{pf}
Let $r(\sF,\Delta)$ be the right hand side of the equality. 
It is clear that $\textnormal{mld}(\sF,\Delta)\leq r(\sF,\Delta)$. 

Let $E$ be the exceptional divisor for some birational morphism $f: Y \rw X$. 
It is well-known that $f$ is a composition of $t$ blowups for some $t\in\bN$. 
Without loss of generality, we assume that $t$ is minimal, which will be denoted as $t(E)$.  
Then it suffices to show the following claim. 
\begin{claim}
For any divisor $E$ over $X$, if $a(E,\sF,\D)+\ve(E)>0$, then $a(E,\sF,\D)+\ve(E) \geq r(\sF,\D)$. 
\end{claim}
\begin{pf}[Claim]
Note that, for $t(E)=0$, we have that $E$ is a divisor on $X$ and 
\[a(E,\sF,\D)+\ve(E) = -\textnormal{ord}_E\D+\ve(E)\geq r(\sF,\D).\] 
When $t(E)=1$, we have $v(E)=0$ and $a(E,\sF,\D)\geq r(\sF,\D)$. 
Then we will proceed by induction on $t$. 

Suppose $t(E)\geq 2$. Nothe that $\ve(E)=0$. 
Let $g_1 : X_1 \rw X$ be the blowup at $f(E)$ and $E_1$ be the (invariant) exceptional divisor for $g_1$.  
Let $K_{\sF_1}+\Delta_1\equiv g_1^*(K_\sF+\Delta)$ where $\sF_1$ is the pullback foliation on $X_1$. 
Then we have 
\[r(\sF_1,\Delta_1) \geq \min\{r(\sF,\Delta), a(E_1,\sF,\Delta)\} \geq r(\sF,\Delta)\]
where the last inequality follows since $t(E_1)=1$. 
Moreover, we have 
\[a(E,\sF,\Delta) = a(E,\sF_1,\D_1) \geq r(\sF_1,\Delta_1)\geq r(\sF,\Delta)\]
where the first inequality comes from the induction hypothesis. 
\end{pf}
\qed
\end{pf}

\begin{cor}\label{mld_min}
Given $(X,\sF,\D= \sum a_iD_i)$ a foliated triple with only log canonical foliation singularities where $a_i\leq 1$. 
Then there is a (log) resolution $f : (Y,\sG,\T) \rw (X,\sF,\Delta)$ such that 
\begin{enumerate}
\item $K_{\sG}+\T = f^*(K_\sF+\Delta)$, 
\item $Y$ is smooth, 
\item $\sG$ has only reduced singularities, 
\item $\T$ has the simple normal crossing support, 
\item the proper transform $\D_Y$ on $Y$ of $\D$ is smooth, 
\item the union of the support of $\T$ and the exceptional divisor of $f$ is $\cup_iT_i$, 
\item all non-invariant irreducible components of $\cup_iT_i$ are pairwisely disjoint, and 
\item any separatrix $C$ through a non-smooth foliation point on a non-invariant prime divisor $T_i$ has $C\cdot T_i=1$. 
\end{enumerate}
Then we have 
{\everymath={\displaystyle}
\[\textnormal{mld}(\sF,\D)=\min\left\{\begin{array}{l}
\min_{E_i\in\wt{I}}\{1+a(E_i,\sF,\Delta)\},\min_{E_i\not\in\wt{I}}\{a(E_i,\sF,\Delta)\}, \\
\min_{(f_*)^{-1}D_i\in\wt{I}}\{1-a_i\},\min_{(f_*)^{-1}D_i\not\in\wt{I}}\{-a_i\}, 1
\end{array}\right\}\]
}
where $\wt{I} = \{T_k\,\vert\, T_k \mbox{ in non-invariant and has only smooth foliation points}\}$. 
\end{cor}
\begin{pf}
First we take $\pi : (Y,\sF_Y) \rw (X,\sF)$ to be the minimal resolution of $(X,\sF)$. 
Then we take $\phi : Z \rw Y$ to be a log resolution of $(Y,(\pi_*)^{-1}\Delta+\textnormal{Exc}(\pi))$. 
After some further blowups, this gives the existence of such $f$. 

Let $b_j = -a(E_j,\sF,\Delta)\leq \ve(E_j)$. 
By Proposition~\ref{mld_sm_model}, we have that 
$\textnormal{mld}(\sG,\T)$ is a minimum of numbers of the following forms:
\begin{itemize}
\item $1$.
\item $-a_i$ if $(f_*)^{-1}D_i\not\in \wt{I}$.
\item $1-a_i$ if $(f_*)^{-1}D_i\in \wt{I}$.
\item $-b_j$ if $E_j\not\in \wt{I}$.
\item $1-b_j$ if $E_j\in \wt{I}$.
\item $-a_i-b_j$ where $(f_*)^{-1}D_i\cap E_j\neq\emptyset$ has some non-smooth foliation singularities. 
\item $1-a_i-b_j$ where $(f_*)^{-1}D_i\cap E_j\neq\emptyset$ has only smooth foliation singularities. 
\item $-b_{j_1}-b_{j_2}$ for $E_{j_1}\cap E_{j_2}\neq\emptyset$ has some non-smooth foliation singularities. 
\item $1-b_{j_1}-b_{j_2}$ for $E_{j_1}\cap E_{j_2}\neq\emptyset$ has only smooth foliation singularities. 
\end{itemize}

Then it suffices to show the following claim. 
\begin{claim}
The last four terms above are irrelevant when taking the minimum. 
\end{claim}
\begin{pf}[Claim]
It is clear if one of $a_i$ and $b_j$ (resp. $b_{j_1}$ and $b_{j_2}$) is non-positive. 
So we may assume that both numbers are positive. 
Moreover, the associated divisors are contained in $\wt{I}$; otherwise, either $-a_i$ or $-b_j$ (resp. either $-b_{j_1}$ or $-b_{j_2}$) is strictly less than zero, which gives a contradiction. 
Thus, both associated divisors are non-invariant, and hence they have empty intersections. 
This completes the proof of the claim.
\end{pf}
\qed
\end{pf}

\section{Ascending chain condition for foliated minimal log discrepancy}
In this section, we fix the following notations. 

The set $B\subset [0,1]$ is always assumed to satisfy the descending chain condition. 
Let $(\sF,\D,x)$ be the germ of a log canonical foliation singularity where $\D$ is a $\bQ$-divisor whose coefficients are in $B$. 
Let $\Gamma=\{E_j\}_{j=1}^r$ be the dual graph of exceptional divisors of a resolution $\pi: (Y,\sG) \rw (X,\sF,x)$ and $\T = \sum_{i=1}^s b_iB_i$ be the proper transform of $\D$ on $Y$. 
Let $a_j := a(E_j,\sF,\D)\geq 0$ all $j=1$, $\ldots$, $r$. 

\begin{lem}\label{basic}
Assume that $\G$ is a tree, all $E_j$ are invariant, and $a_j\leq 1$ for all $j=1$, $\ldots$, $r$. 
Then we have:
\begin{enumerate}[(a)]
\item If $a_j \geq \ve$ for some positive real number $\ve$, then $-E_j^2 \leq\lfloor\frac{2}{\ve}\rfloor$. 
\item If $-E_j^2\geq 2$ for some $j$, then 
\[1-a_j \geq \frac{1-a_k}{2}\] 
and $2a_j\leq a_k + a_\ell$ for any $k\neq \ell$ such that $E_j\cdot E_k = E_j\cdot E_\ell = 1$. 
\item Suppose the vertex $[E_{j_0}]$ is a fork and the vertices $[E_{j_k}]$ are connected to $[E_{j_0}]$ for $k=1$, $2$, $3$. 
If $-E_{j_k}^2\geq 2$ for $k=0$, $1$, $2$, then $a_{j_0}\leq a_{j_3}$. 
\item If $-E_j^2=1$ for some $j$, $-E_k^2\geq 2$, $a_j<a_k$ for all $k\neq j$, and the vertex $[E_j]$ is not a fork, then $\Gamma$ is a chain. 
\item Given a sequence of vertices $[E_1]$, $\ldots$, $[E_m]$ where $E_i\cdot E_{i+1}=1$ for all $i=1$, $\ldots$, $m-1$. 
Suppose $-E_i^2\geq 2$ for $i = 2$, $\ldots$, $m-1$ and $a_1\leq a_2$. 
If either $-E_2^2\geq 3$ and $a_2\geq \ve$ or $\T\cdot E_2\geq \ve$, then $m\leq\lfloor\frac{1}{\ve}\rfloor + 2$. 
\end{enumerate}
\end{lem}

\begin{pf}
Notice that 
\begin{equation}\label{basic_ineq}
\left\{\begin{array}{l}
\begin{aligned}
0 &= \left(K_\sG+\T - \sum_{i=1}^r a_iE_i\right)\cdot E_j \\
&= K_\sG\cdot E_j + \sum_{i=1}^s t_{i,j}b_i - a_jE_j^2 - \sum_{E_i\cdot E_j=1} a_i \\
&\geq  \#\{i\,\vert\, E_i\cdot E_j=1\} -2 + \sum_{i=1}^s t_{i,j}b_i - a_jE_j^2 - \sum_{E_i\cdot E_j=1} a_i \\
&= -2 + \sum_{i=1}^s t_{i,j}b_i - a_jE_j^2 + \sum_{E_i\cdot E_j=1}(1-a_i)
\end{aligned}
\end{array}\right.
\end{equation}
where $t_{i,j} =B_i\cdot E_j\geq 0$.
\begin{enumerate}[(a)]
\item From the inequality~(\ref{basic_ineq}), we have that $0\geq -2 - a_jE_j^2$. 
Thus, \[-E_j^2 \leq \frac{2}{a_j} \leq\frac{2}{\ve}.\]
Since $-E_j^2\in\bZ$, we have that $-E_j^2 \leq\lfloor\frac{2}{\ve}\rfloor$. 
\item From the inequality~(\ref{basic_ineq}), we have that $0\geq -2 - a_jE_j^2 + 1 - a_k\geq -1+2a_j-a_k$ and $0\geq -2 - a_jE_j^2 + 1 - a_k+1-a_\ell\geq 2a_j-a_k-a_\ell$. 
Then we get $1-a_j \geq \frac{1-a_k}{2}$ and $2a_j\leq a_k + a_\ell$. 
\item By (b), we have 
\[1-a_{j_k}\geq \frac{1-a_{j_0}}{2}\]
for $k=1$, $2$. Thus 
\begin{eqnarray}\label{fork}
1-a_{j_1}+1-a_{j_2} \geq 1-a_{j_0}. 
\end{eqnarray}
Also from the inequality~(\ref{basic_ineq}), we have 
\[0\geq 1 +2a_{j_0} - \sum_{k=1}^3 a_{j_k}. \]
Combining with the inequality~(\ref{fork}), we get 
\[a_{j_3}\geq 1+2a_{j_0}-a_{j_1}-a_{j_2}\geq 2a_{j_0}-a_{j_0} = a_{j_0}.\]
\item Assume $\G$ is not a chain. 
Let the vertex $[E_k]$ be a fork. 
Let $[E_k] =: [E_{\ell_0}]$, $[E_{\ell_1}]$, $\ldots$, $[E_{\ell_m}] := [E_j]$ be a sequence of vertex connecting from $[E_k]$ to $[E_j]$. 
By (c), we know that $a_{\ell_0}\leq a_{\ell_1}$. 
By (b), we get $a_{\ell_1}\leq a_{\ell_2} \leq \ldots\leq a_{\ell_m}$. 
Thus, $a_k\leq a_j$, which gives a contradiction. 
\item From the inequality~(\ref{basic_ineq}) with $j=2$, we have that 
\begin{eqnarray*}
0 &\geq & -2+\sum_{i=1}^s t_{i,2}b_i-a_2E_2^2 + (1-a_1)+(1-a_3) \\
&\geq & 2a_2-a_1-a_3+\ve \\
&\geq & a_2-a_3+\ve
\end{eqnarray*}
Therefore, $a_3\geq a_2+\ve$. 
\begin{claim}
$a_{j+1}-a_j \geq \ve$ for all $j = 2$, $\ldots$, $m-1$. 
\end{claim}
\begin{pf}[Claim]
We have seen the claim holds true when $j=2$. 
Then by (b), we have $a_{j+1}+a_{j-1}\geq 2a_j$. 
Thus, $a_{j+1}-a_j\geq a_j-a_{j-1}\geq \ve$ by induction on $j$.  
\end{pf}
Therefore, we have 
\[1\geq a_m \geq a_{m-1}+\ve \geq \cdots\geq a_2+(m-2)\ve \geq (m-2)\ve.\]
Hence, $m\leq\lfloor\frac{1}{\ve}\rfloor + 2$. 
\end{enumerate}
\qed
\end{pf}

\begin{lem}\label{lem_finite}
Assume that $(Y,\sG)$ is the minimal resolution of the germ $(X,\sF,x)$. Then we have the followings:
\begin{enumerate}[(a)]
\item $\textnormal{pld}_x(\sF) = 0$ if $\G$ is not of the type (1) in Theorem~\ref{main}. 
\item If $a_j>0$ for all $j$, then $a_j \leq a(E_j,\sF)\leq 1$ for all $E_j$'s. 
\item Fix a number $\delta$, there are only finitely many sequences $\{t_1$, $\ldots$, $t_r\}$ where $t_i\in\bN$ such that $\sum_{i=1}^r t_ib_i \leq \delta$ for some $b_i\in B\backslash\{0\}$. 
\end{enumerate}
\end{lem}
\begin{pf}
\begin{enumerate}[(a)]
\item This is straightforward. 
\item Note that 
\[\left(\sum_{i=1}^r a_iE_i\right)\cdot E_j = (K_\sG+\T)\cdot E_j\geq K_\sG\cdot E_j = \left(\sum_{i=1}^r a(E_i,\sF)E_i\right)\cdot E_j\]
for all $j$. 
By Lemma~\ref{key_lem}, we have $0<a_j \leq a(E_j,\sF)$ for all $j$. 
Then $\G$ is of type (1) in Theorem~\ref{main}. 
Also, we notice that 
\[K_\sG\cdot E_j\geq \left(\sum_{i=1}^r E_i\right)\cdot E_j.\] 
Thus, we have $a(E_j,\sF)\leq 1$ by Lemma~\ref{key_lem}.

\item Since $B$ is a DCC set, there is a positive number $\ve$ such that $b_i\geq \ve$ for all $b_i\in B\backslash\{0\}$. 
Note that 
\[\delta \geq \sum_{i=1}^r t_ib_i \geq \sum_{i=1}^r t_i\ve \mbox{ and thus } \sum_{i=1}^r t_i\leq \frac{\delta}{\ve}.\]
This shows that there are only finitely many possible $r$ and, for any fixed $r$, there are only finitely many sequences $\{t_1$, $\ldots$, $t_r\}$ such that $\sum t_i\leq \frac{\delta}{\ve}$. 
This proves (c).
\end{enumerate}
\qed
\end{pf}

\begin{lem}\label{pld_classification}
Fix $\ve>0$. 
Suppose $(Y,\sG)$ is the minimal resolution of the germ $(X,\sF,x)$ with $\textnormal{pld}_x(\sF,\D)\geq \ve$ and $b_j\geq\ve$ for all $j$. 
Then $\G$ belongs to one of the following cases:
\begin{enumerate}
\item Finitely many graphs (that include the way how $B_i$ intersects $E_j$).
\item The chain $\cup_jE_j$ given by the ordered curves $L_{\ell_1}$, $\ldots$, $L_1$, $M_1$, $\ldots$, $M_n$, $R_1$, $\ldots$, $R_{\ell_2}$ where $K_\sG\cdot L_{\ell_1}=-1$, the weights of $M_k$'s are 2, and each $B_j$ does not meet any $M_k$'s. 
Moreover, the partial log discrepancy $\textnormal{pld}_x(\sF,\D)$ is achieved at either $L_1$ or $R_1$ and  there are only finitely many possibilities (independent of $n$) for the dual graphs $\{L_1$, $\ldots$, $L_{\ell_1}\}$ and $\{R_1$, $\ldots$, $R_{\ell_2}\}$ and the way how $B_i$ intersects $L_\alpha$ and $R_\beta$. 
\end{enumerate}
\end{lem}

\begin{pf}
By Lemma~\ref{lem_finite},  we have that $\G$ is of type (1), and therefore $\G$ satisfies the assumption in Lemma~\ref{basic}. 
So from the inequality~(\ref{basic_ineq}), we have that $0\geq -2 + \sum t_{i,j}b_i -a_jE_j^2$ and hence $\sum t_{i,j}b_i\leq 2+a_jE_j^2\leq 2$.
Therefore, there are only finitely many possibilities for $t_{i,j}$ for any $i$, $j$ by Lemma~\ref{lem_finite} (c). 

Let the chain be $\cup_jE_j$. 
By Lemma~\ref{basic} (b), the function $f: \{1$, $\ldots$, $r\} \rw \bR$ which maps $j$ to $a_j$ is convex. 
Let $S = \{j\vert\, a_j=\min_k\{a_k\}\}$, $N=\#S$, and $j_0=\min S$.  

If $N\geq 2$, then $S = \{j_0$, $\ldots$, $j_0+N-1\}$ by convexity of $f$.
Note that, by inequality~(\ref{basic_ineq}) with fixed $j$, we have 
\[0\geq -2 + \sum_{i=1}^s t_{i,j}b_i - a_jE_j^2 + (1-a_{j-1})+(1-a_{j+1})\]
and thus 
\begin{equation}\label{eq_conv}
a_{j-1}-2a_j+a_{j+1}\geq \sum_{i=1}^s t_{i,j}b_i-a_j(E_j^2+2)\geq\left(\sum_{i=1}^s t_{i,j}-E_j^2-2\right)\ve\geq \ve
\end{equation}
whenever $-E_j^2\geq 3$ or $t_{i,j}\geq 1$ for some $i$. 
So we can choose $L_1 = E_{j_0}$ and $R_1 = E_{j_0+N-1}$. 
By Lemma~\ref{basic} (e), we have $\ell_i\leq \lfloor\frac{2}{\ve}\rfloor+1$ for $i=1$, $2$. 

If $N=1$, then $S = \{j_0\}$. 
Let $T = \{j\vert\, E_j^2=-2 \mbox{ and } E_j\cdot\T=0\}$. 
We have the following three cases:
\begin{enumerate}[(i)]
\item If both $j_0-1$ and $j_0+1$ are not in $T$, then we choose $L_1 = E_{j_0}$ and $R_1 = E_{j_0+1}$. 
\item If $j_0-1\in T$, then we choose $R_1 = E_{j_0}$ and $L_1 = E_u$ where $u$ is the maximal integer strictly less than $j_0$ and not in $T$. 
\item If $j_0+1\in T$, then we choose $L_1 = E_{j_0}$ and $R_1 = E_u$ where $u$ is the minimal integer strictly greater than $j_0$ and not in $T$. 
\end{enumerate}
By Lemma~\ref{basic} (e), we also have $\ell_i\leq \lfloor\frac{2}{\ve}\rfloor+1$ for $i=1$, $2$. 
To complete the proof, it remains to show the following claim. 

\begin{claim}
There are only finitely many possibilities for the dual graphs $\{L_1$, $\ldots$, $L_{\ell_1}\}$ and $\{R_1$, $\ldots$, $R_{\ell_2}\}$ and the way how $B_i$ intersects $L_\alpha$ and $R_\beta$. 
\end{claim}
\begin{pf}
To simplify some notations, we assume that we have a chain $C = \cup_{j=1}^\ell C_j$ with increasing associated discrepancies $a_j\geq\ve$. 
By inequality~(\ref{eq_conv}), we have 
\begin{eqnarray*}
1 &\geq & a_\ell \\
&\geq & a_{\ell-1} + \left(\sum_{i=1}^s t_{i,\ell-1}-E_{\ell-1}^2-2\right)\ve \\
&\vdots & \\
&\geq & a_1 + \sum_{j=1}^{\ell-1}\left(\sum_{i=1}^s t_{i,j}-E_{j}^2-2\right)\ve \\
&\geq & \ve + \sum_{j=1}^{\ell-1}\left(\sum_{i=1}^s t_{i,j}-E_{j}^2-2\right)\ve.
\end{eqnarray*}
Thus, we have \[\sum_{j=1}^{\ell-1}\left(\sum_{i=1}^s t_{i,j}-E_{j}^2-2\right)\leq \frac{1}{\ve}-1.\]
Therefore, there are only finitely many possibilities for $E_j^2$ for any $j$ since we have seen that there are only finitely many possibilities for $t_{i,j}$ for any $i$, $j$.
\end{pf}
\qed
\end{pf}

\begin{lem}\label{pld}
For the second case in Lemma~\ref{pld_classification}, if the chain has length $r = \ell_1+n+\ell_2 > 2n_0+2$ 
where $n_0 = \lfloor \frac{1}{\ve}\rfloor$, 
then \[\textnormal{mld}_x(\sF,\Delta) = \textnormal{pld}_x(\sF,\Delta).\]

Moreover, for the fixed graphs $\{L_1$, $\ldots$, $L_{\ell_1}\}$ and $\{R_1$, $\ldots$, $R_{\ell_2}\}$, the fixed number $t$ of irreducible components of $B$, and the fixed way how $B_i$ intersects $E_j$, if $\lim_{n\to\infty}\, b_j = \overline{b_j}$ exists where $(b_1$, $\ldots$, $b_s)$ is some ordering of coefficients of $\T$, then 
\[\textnormal{pld}(\sF,\Delta)\geq\min\left\{\frac{\alpha^L}{m_1-q_1},\frac{\alpha^R}{m_2-q_2}\right\}\] and 
\[\lim_{n\to\infty}\,\textnormal{pld}(\sF,\Delta) = \min\left\{\frac{\overline{\alpha^L}}{m_1-q_1}, \frac{\overline{\alpha^R}}{m_2-q_2}\right\}.\] 
where 
\begin{itemize}
\item $c_i^L = \sum_j(B_i\cdot L_j)g^L_j$ and $c_i^R = \sum_j(B_i\cdot R_j)g^R_j$, 
\item $g^L_j$ (resp. $g^R_j$) is the determinant of the chain $L_{j+1}$, $\ldots$, $L_{\ell_1}$ 
(resp. $R_{j+1}$, $\ldots$, $R_{\ell_2}$), 
\item $m_1 = g^L_0$, $q_1 = g^L_1$, $m_2 = g^R_0$, and $q_2 = g^R_1$,
\item $\alpha^L = 1-\sum_ib_ic_i^L$ and $\alpha^R = -\sum_ib_ic_i^R$,
\item $\overline{\alpha^L} = 1-\sum_i\overline{b_i}c_i^L$ and $\overline{\alpha^R} = -\sum_i\overline{b_i}c_i^R$. 
\end{itemize}

\end{lem}
\begin{pf}
We first show that $\textnormal{mld}_x(\sF,\Delta) = \textnormal{pld}_x(\sF,\Delta)$ if $r> 2n_0+2$. 
Suppose not, then there is an irreducible divisor $F$ over $Y$ such that $a(F,\sF,\Delta)<\textnormal{pld}(\sF,\Delta)$. 
Then $F$ is an exceptional divisor for the composition of $t$ blowups. 
Let $F_1$, $\ldots$, $F_t:=F$ be the exceptional divisors for the corresponding blowups.
We may assume that $a(F,\sF,\Delta) < a(F_j,\sF,\Delta)$ for all $1\leq j\leq t-1$. 
Note that $-F_j^2\geq 2$ for all $1\leq j\leq t-1$ and $[F]$ is not a fork since $\sum_{i=1}^{t-1}F_i$ has the simple normal crossing support. 
Then by Lemma~\ref{basic} (d), the graph $\{E_i,F_j\} = \{G_k\}_{k=1}^N$ is a chain with $F=F_t=G_u$ for some $u$. 
Notice that exactly one of $-G^2_{u-1}$ or $-G^2_{u+1}$ is 2. 
Without loss of generality, we assume that $-G^2_{u+1} = 2$ and $-G^2_{u-1}\geq 3$. 
Then $u\leq n_0+2$ by Lemma~\ref{basic} (e). 

Put $\gamma = \min\{k> u\vert -G^2_k\geq 3\}$. 
Then, by Lemma~\ref{basic} (e) again, we have $N-\gamma\leq n_0$. 
By blowing down $G_u$, $\ldots$, $G_{\gamma-1}$ successively, we get a new foliated surface with only reduced singularities and with at most $(u-1)+(N-\gamma+1)\leq 2n_0+2<r$ exceptional divisors over $X$. 
However, since this new foliated surface factors through the minimal resolution, the number of exceptional divisors over $X$ is at least $r$, which gives a contradiction. 

Now we compute the partial log discrepancy $\textnormal{pld}_x(\sF,\Delta)$. 
It is known that the graph $\G$ is uniquely determined by these five numbers $m_1$, $q_1$, $m_2$, $q_2$, $n$. 
By Lemma~\ref{pld_calc}, we have that both $a(L_1,\sF,\D)$ and $a(R_1,\sF,\D)$ are between 
\[\min\left\{\frac{\alpha^L}{m_1-q_1},\frac{\alpha^R}{m_2-q_2}\right\} \mbox{ and } \max\left\{\frac{\alpha^L}{m_1-q_1},\frac{\alpha^R}{m_2-q_2}\right\}.\] 
Moreover, we have that 
\[\lim_{n\to\infty}\,a(L_1,\sF,\D) = \frac{\overline{\alpha^L}}{m_1-q_1} \mbox{ and } \lim_{n\to\infty}\,a(R_1,\sF,\D) = \frac{\overline{\alpha^R}}{m_2-q_2}.\]
\qed
\end{pf}

\begin{lem}\label{pld_calc}
Notation as in Lemma~\ref{pld}, we have 
\begin{eqnarray*}
a(L_1,\sF,\Delta) &=& 
\frac{\alpha^L(n(m_2-q_2)+m_2) + \alpha^Rq_1}{n(m_1-q_1)(m_2-q_2) + m_2(m_1-q_1) + q_1(m_2-q_2)} \\
&=& \frac{\frac{\alpha^L}{m_1-q_1}(n+\frac{m_2}{m_2-q_2})+\frac{\alpha^R}{m_2-q_2}\frac{q_1}{m_1-q_1}}{n+\frac{m_2}{m_2-q_2}+\frac{q_1}{m_1-q_1}}
\end{eqnarray*}
and 
\begin{eqnarray*}
a(R_1,\sF,\Delta) &=& 
\frac{\alpha^R(n(m_1-q_1)+m_1) + \alpha^Lq_2}{n(m_1-q_1)(m_2-q_2) + m_1(m_2-q_2) + q_2(m_1-q_1)} \\
&=& \frac{\frac{\alpha^R}{m_2-q_2}(n+\frac{m_1}{m_1-q_1})+\frac{\alpha^L}{m_1-q_1}\frac{q_2}{m_2-q_2}}{n+\frac{m_1}{m_1-q_1}+\frac{q_2}{m_2-q_2}}.
\end{eqnarray*}
\end{lem}
\begin{pf}
We denote $\G$ as $\G_{m_1,q_1,m_2,q_2,n}$ and its intersection matrix $A$ as $A_{m_1,q_1,m_2,q_2,n}$. 
Let $S_n = \det(-A_{1,1,m_2,q_2,n})$. 
Then by Lemma~\ref{det} (1), we have that 
\[\det(-A_{m_1,q_1,m_2,q_2,n}) = 2m_1S_{n-1} - a_1S_{n-1} - m_1S_{n-2} = (2m_1 - q_1)S_{n-1} - m_1S_{n-2}\] 
and $S_n = 2S_{n-1}-S_{n-2}$. 
Thus, $S_{i+1}-S_i = S_{i}-S_{i-1}$ for all $i$. 
Since $S_0 = m_2$ and $S_1 = 2S_0 - q_2 = 2m_2 - q_2$, we have 
\[S_n = \sum_{i=1}^n(S_i-S_{i-1}) + S_0 = n(S_1-S_0) + S_0 = n(m_2-q_2)+m_2\] 
and thus 
\begin{eqnarray*}
\det(-A) &=& \det(-A_{m_1,q_1,m_2,q_2,n}) \\
&=& (2m_1-q_1)\big((n-1)(m_2-q_2)+m_2\big) - m_1\big((n-2)(m_2-q_2)+m_2\big) \\
&=& n(m_1-q_1)(m_2-q_2) + m_2(m_1-q_1) + q_1(m_2-q_2).
\end{eqnarray*}

Also by Lemma~\ref{det} (2), we have 
\begin{eqnarray*}
(L_1,L_j)\textnormal{-cofactor} &=& (-1)^{|\G|+1} g^L_jS_n, \\ 
(L_1,M_j)\textnormal{-cofactor} &=& (-1)^{|\G|+1} q_1S_{n-j}, \mbox{ and} \\
(L_1,R_j)\textnormal{-cofactor} &=& (-1)^{|\G|+1} g^R_jq_1.
\end{eqnarray*}
Put $d_j^L = (K_\sG+\T)\cdot L_j$, $d_j^R = (K_\sG+\T)\cdot R_j$, and $d_j^M = (K_\sG+\T)\cdot M_j$. 
By assumption, we know that $d_j^M=0$ and $d_j^R = \T\cdot R_j$ for all $j$. 
Moreover, $d_j^L = \T\cdot L_j$ for $j\neq\ell_1$ and $d_{\ell_1}^L = -1+\T\cdot L_{\ell_1}$. 
Then we have 
\[\det(A)a(L_1,\sF,\Delta) = \sum_{j=1}^{\ell_1}d_j^L\cdot(L_1,L_j)\textnormal{-cofactor} +\sum_{j=1}^{\ell_2}d_j^R\cdot(L_1,R_j)\textnormal{-cofactor}. \]
Thus 
\begin{eqnarray*}
& & -\det(-A)a(L_1,\sF,\Delta) \\
&=& \sum_{j=1}^{\ell_1}d_j^L g^L_jS_n + \sum_{j=1}^{\ell_2}d_j^Rg^R_jq_1 \\
&=& \sum_{j=1}^{\ell_1}\left(\sum_{i=1}^sb_i(B_i\cdot L_j)\right)g^L_jS_n - g^L_{\ell_1}S_n + \sum_{j=1}^{\ell_2}\left(\sum_{i=1}^sb_i(B_i\cdot R_j)\right)g^R_jq_1 \\
&=& \left(\sum_{i=1}^sb_ic_i^L-1\right)S_n + \left(\sum_{i=1}^sb_ic_i^R\right)q_1\\
&=& -\alpha^LS_n-\alpha^Rq_1.
\end{eqnarray*}
Therefore, 
\begin{eqnarray*}
a(L_1,\sF,\Delta) &=& 
\frac{\alpha^L\big(n(m_2-q_2)+m_2\big) + \alpha^Rq_1}{n(m_1-q_1)(m_2-q_2) + m_2(m_1-q_1) + q_1(m_2-q_2)} \\
&=& \frac{\frac{\alpha^L}{m_1-q_1}(n+\frac{m_2}{m_2-q_2})+\frac{\alpha^R}{m_2-q_2}\frac{q_1}{m_1-q_1}}{n+\frac{m_2}{m_2-q_2}+\frac{q_1}{m_1-q_1}}.
\end{eqnarray*}
By the similar way, we also have 
\begin{eqnarray*}
a(R_1,\sF,\Delta) &=& 
\frac{\alpha^R(n(m_1-q_1)+m_1) + \alpha^Lq_2}{n(m_1-q_1)(m_2-q_2) + m_1(m_2-q_2) + q_2(m_1-q_1)} \\
&=& \frac{\frac{\alpha^R}{m_2-q_2}(n+\frac{m_1}{m_1-q_1})+\frac{\alpha^L}{m_1-q_1}\frac{q_2}{m_2-q_2}}{n+\frac{m_1}{m_1-q_1}+\frac{q_2}{m_2-q_2}}.
\end{eqnarray*}
\qed
\end{pf}

\begin{thm}\label{pld_acc}
For any DCC set $B$, the set 
\[\textnormal{PLD}(2,B) := \{\textnormal{pld}_x(\sF,\Delta) \vert \, \mbox{ $(X,\sF,\Delta)$ is a foliated triple with $x\in X$ and $\Delta\in B$}\}\] 
satisfies the ascending chain condition (ACC). 
\end{thm}
\begin{pf}
Given any non-decreasing sequence $\{\textnormal{pld}_{x_k}(\sF_k,\Delta_k)\}_{k=1}^\infty$ in the set $\textnormal{PLD}(2,B)$ where $(X_k,\sF_k,\Delta_k)$ is a germ of foliated triple around $x_k$ and $\D_k \in B$ for all $k$. 
We may assume that 
\[\textnormal{pld}_{x_k}(\sF_k,\Delta_k) > 0\] 
for all $k$, otherwise the sequence $\{\textnormal{pld}_{x_k}(\sF_k,\Delta_k)\}_{k=1}^\infty$ stabilizes. 
Now let $\ve>0$ be a number such that $\textnormal{pld}_{x_k}(\sF_k,\Delta_k) \geq \ve$ for all $k$. 
Since $B$ satisfies the descending chain condition, we may assume that $\min(B\backslash\{0\})\geq\ve$. 
\begin{claim}
We may also assume that the number $s$ of irreducible components of $\T_k$ is independent of $k$. 
\end{claim}
\begin{pf}
For ease of our notation, we will drop the subscription $k$. 

By Lemma~\ref{pld_classification}, the number of irreducible components of exceptional divisors $E_i$ which meet $\T$ is bounded by $N$ for some $N>0$. 
Moreover, by the inequality~(\ref{basic_ineq}), we have $\sum_{i=1}^s t_{i,j}b_i\leq 2$ where $s$ is the number of irreducible components of $\T$. Then 
\[2N\geq \sum_{j=1}^r\sum_{i=1}^s t_{i,j}b_i \geq \sum_{i=1}^s\sum_{j=1}^r t_{i,j}\ve\geq\sum_{i=1}^s\ve = s\ve\]
where the third inequality follows from $t_{i,j}\geq 1$ for some $j$ for any fixed $i$. 
Therefore, $s\leq\frac{2N}{\ve}$. 
Hence, after taking a subsequence, we may assume that the number of irreducible components of $\T_k$ is independent of $k$. 
\end{pf}

Let $(b_{1,k}$, $\ldots$, $b_{s,k})$ be some ordering of coefficients of $\T_k$. 
After taking a further subsequence, we may assume the sequence $\{b_{j,k}\}_{k=1}^\infty$ is non-decreasing for any $j = 1$, $\ldots$, $s$. 
If there are infinitely many dual graph $\G_k$ belonging to the first case of Lemma~\ref{pld_classification}, then, by Lemma~\ref{key_lem} and after taking a further subsequence, we have that $\textnormal{pld}_{x_k}(\sF_k,\Delta_k)$ is non-increasing. 
Since it is also non-decreasing, the sequence $\{\textnormal{pld}_{x_k}(\sF_k,\Delta_k)\}_{k=1}^\infty$ stabilizes. 

Thus, there are infinitely many dual graphs $\G_k$ belonging to the second case of Lemma~\ref{pld_classification}. 
After taking a subsequence, we may assume that the dual graphs of $\{L_1$, $\ldots$, $L_{\ell_1}\}$ and $\{R_1$, $\ldots$, $R_{\ell_2}\}$ are fixed, the way how $B_i$ intersects $E_j$ is fixed, and $n$ is sufficiently large. 
Since $b_{j,k}\leq 1$, the limit $\lim_{k\to\infty}\, b_{j,k} = \overline{b_j}$ exists for any $j = 1$, $\ldots$, $s$. 

By Lemma~\ref{pld}, we have, for all $k$, that 
\[\textnormal{pld}_{x_k}(\sF_k,\Delta_k)\geq \min\left\{\frac{\alpha^L}{m_1-q_1},\frac{\alpha^R}{m_2-q_2}\right\}
\geq\min\left\{\frac{\overline{\alpha^L}}{m_1-q_1}, \frac{\overline{\alpha^R}}{m_2-q_2}\right\}\]
and \[\lim_{k\to\infty}\,\textnormal{pld}_{x_k}(\sF_k,\Delta_k) 
= \min\left\{\frac{\overline{\alpha^L}}{m_1-q_1}, \frac{\overline{\alpha^R}}{m_2-q_2}\right\}.\]

Thus, the sequence $\{\textnormal{pld}_{x_k}(\sF_k,\Delta_k)\}_{k=1}^\infty$ stabilizes. 
\qed
\end{pf}

\begin{thm}\label{mld_acc}
For any DCC set $B$, the set 
\[\textnormal{MLD}(2,B) := \{\textnormal{mld}_x(\sF,\Delta) \vert \, \mbox{ $(X,\sF,\Delta)$ is a foliated triple with $x\in X$ and $\Delta\in B$.}\}\] 
satisfies the ascending chain condition (ACC). 
\end{thm}
\begin{pf}
Given any non-decreasing sequence $\{\textnormal{mld}_{x_k}(\sF_k,\Delta_k)\}_{k=1}^\infty$ in the set $\textnormal{MLD}(2,B)$ where $(X_k,\sF_k,\Delta_k)$ is a germ of foliated triple around $x_k$ and $\D_k \in B$ for all $k$. 
We may assume that 
\[\textnormal{pld}_{x_k}(\sF_k,\Delta_k) > 0\] 
for all $k$, otherwise the sequence 
$\{\textnormal{mld}_{x_k}(\sF_k,\Delta_k)\}_{k=1}^\infty$ stabilizes. 
Now let $\ve>0$ be a number such that $\textnormal{mld}_{x_k}(\sF_k,\Delta_k) \geq \ve$ for all $k$. 
Since $B$ satisfies the descending chain condition, we may assume that $\min(B\backslash\{0\})\geq\ve$. 

If there are infinitely many $k$ such that $\textnormal{pld}_{x_k}(\sF_k,\Delta_k)=\textnormal{mld}_{x_k}(\sF_k,\Delta_k)$, then by Theorem~\ref{pld_acc}, we get the sequence $\{\textnormal{mld}_{x_k}(\sF_k,\Delta_k)\}_{k=1}^\infty$ stabilizes. 

Thus, by Lemma~\ref{pld_classification} and \ref{pld}, after taking a subsequence, we may assume $(X_k,\sF_k,\Delta_k)$ have the same weighted dual graph for the minimal resolution $(Y_k,\sG_k,\T_k)$, the ways how $B_i$ intersects $E_j$ are the same, and $\textnormal{pld}_{x_k}(\sF_k,\Delta_k)>\textnormal{mld}_{x_k}(\sF_k,\Delta_k)$. 

We may assume the number of irreducible components of $\T_k$ is fixed as the claim in the proof of Theorem~\ref{pld_acc}. 
Let $(b_{1,k}$, $\ldots$, $b_{s,k})$ be some ordering of the coefficients of $\T_k$. 
Then, by taking a further subsequence, we may assume the sequence $\{b_{j,k}\}_{k=1}^\infty$ is non-decreasing for any $j=1$, $\ldots$, $s$. 

Since $\textnormal{pld}_{x_k}(\sF_k,\Delta_k)>\textnormal{mld}_x(\sF_k,\Delta_k)$, there exists an exceptional divisor $F_k$ over $Y_k$ such that $a(F_k,\sF_k,\Delta_k) = \textnormal{mld}_{x_k}(\sF_k,\Delta_k)$. 
Then $F_k$ is an exceptional divisor of the birational morphism $\pi_k: Z_k \rw Y_k$ which is the composition of $N_k$ blowups. 
Let $F_{1,k}$, $\ldots$, $F_{N_k,k} =: F_k$ be all $\pi_k$-exceptional divisors. 
We may assume that $a(F_{j,k},\sF_k,\Delta_k)>a(F_k,\sF_k,\Delta_k)$ for all $j\leq N_k-1$. 
Notice that $-(F_{j,k})^2\geq 2$ for all $j \leq N_k-1$ and $[F_k]$ is not a fork since $\bigcup_iE_i\cup\bigcup_j F_{j,k}$ has the simple normal crossing support. 
By Lemma~\ref{basic} (d), the dual graph of $\{E_i,F_{j,k}\}_{i,j}$ is a chain. 

By Lemma~\ref{bdd_blowup}, the length of the chain of the dual graph of $\{E_i,F_{j,k}\}_{i,j}$ is bounded. 
After taking a subsequence, we may assume they have the same dual graphs, and the ways $B_i$ intersecting $E_j$ and $F_{j,k}$ are the same. 

Therefore, by Lemma~\ref{key_lem}, we have that 
\[\textnormal{mld}_{x_k}(\sF_k,\Delta_k) = a(F_k,\sF_k,\Delta_k)\geq a(F_{k+1},\sF_{k+1},\Delta_{k+1}) 
= \textnormal{mld}_x(\sF_{k+1},\Delta_{k+1})\]
for all $k$. 
This shows that the sequence $\{\textnormal{mld}_{x_k}(\sF_k,\Delta_k)\}_{k=1}^\infty$ stabilizes. 
\qed
\end{pf}

\begin{lem}\label{bdd_blowup}
Fix an $\ve>0$. 
Suppose the sequence $(X_k,\sF_k,\Delta_k)$ with $\textnormal{pld}_{x_k}(\sF_k,\D_k)>\textnormal{mld}_{x_k}(\sF_k,\D_k)\geq\ve$ and $b_i\geq\ve$ have the following properties:
\begin{enumerate}[(i)]
\item The dual graph for the minimal resolution $(Y_k,\sG_k,\T_k)$ are the same.
\item The number $s$ of irreducible components of $\T_k$ are the same. 
\item $(b_{1,k}$, $\ldots$, $b_{s,k})$ is some ordering of the coefficients of $\T_k$. 
\item $\{b_{j,k}\}_{k=1}^\infty$ is non-decreasing for any $j=1$, $\ldots$, $s$. 
\end{enumerate}
Then there is a positive integer $N$ independent of $k$ such that, for each $(Y_k,\sG_k,\T_k)$, there exists a birational morphism $\pi_k : Z_k \rw Y_k$ such that 
\begin{enumerate}
\item the relative Picard number $\rho(Z_k/Y_k)\leq N$ and 
\item one of the exceptional divisors on $Z_k$ over $X_k$ computes the minimal log discrepancy. 
\end{enumerate}
\end{lem}

\begin{pf}
To simplify our notation, we will drop the subscript $k$ and, for any fixed divisor $D$, denote all of the proper transforms of $D$ by $D$. Since the set $B$ satisfies DCC, there is a $\delta_a>0$ for any $a\in\bR$ such that $\sum_{i=1}^sm_i\beta_i-1\geq\delta_a$ if $\sum_{i=1}^sm_i\beta_i-a>0$ for some $\beta_i\in B$. 

We have seen in the proof of Theorem~\ref{mld_acc} that the dual graph of $\{G_n\}=\{E_i,F_j\}$ is a chain. 
This implies that $\pi$ is the composition of blowups with center $p$ either at the foliation singularities or on two curves at the ends. 
We may also assume that $N$ is minimal for all $k$. 
We will proceed induction on $\sum_{i=1}^s\textnormal{mult}_pB_i$. 
Note that 
\[\sum_{i=1}^s\textnormal{mult}_pB_i\leq\frac{1}{\ve}\sum_{j=1}^r\T\cdot E_j\leq \frac{1}{\ve}\sum_{j=1}^r2 = \frac{2r}{\ve}\]
by the inequality~(\ref{basic_ineq}) and $\frac{2r}{\ve}$ is independent of $k$. 

Let $F_1$ be the exceptional divisor of blowup at $p$. 
We have following three cases:
\begin{enumerate}
\item Suppose $p$ is a smooth foliation point on the curve at the end. 
Note that 
\[a(F_1,\sF,\D) = a(E_r,\sF,\D) + 1-\sum_{i=1}^s(\textnormal{mult}_pB_i)b_i.\]

If all $B_i$'s meet $E_r$ at $p$ transversally and $\sum_{i=1}^s(\textnormal{mult}_pB_i)b_i\leq 1$, then the exceptional divisors from further blowups have discrepancies at least $a(F_1,\sF,\D)$.

If all $B_i$'s meet $E_r$ at $p$ transversally but $\sum_{i=1}^s(\textnormal{mult}_pB_i)b_i > 1$, then 
\[\sum_{i=1}^s(\textnormal{mult}_pB_i)b_i -1\geq \delta_1>0\] 
and $a(F_1,\sF,\D) \leq a(E_r,\sF,\D)-\delta_1\leq 1-\delta_1$. 
If $\sum_{i=1}^s\textnormal{mult}_qB_i<\sum_{i=1}^s\textnormal{mult}_pB_i$ for all $q\in F_1$, then we are done by induction. 
If not, then the exceptional divisor from blowing up at such $q$ has discrepancy $\leq a(F_1,\sF,\D)-\delta_1\leq 1-2\delta_1$. 
Since the minimal log discrepancy $\textnormal{mld}_x(\sF,\D)$ is positive, after finitely many blowups, the quantity $\sum_{i=1}^s\textnormal{mult}_pB_i$ will strictly decrease. 

If some of $B_i$'s meet $E_r$ at $p$ transversally and some do not, then the quantity $\sum_{i=1}^s\textnormal{mult}_pB_i$ strictly decreases. 

If each $B_i$ does not meet $E_r$ at $p$ transversally, then all $B_i$'s meet $F_1$ at the reduced singularity $q$ on $F_1$. 
Notice that we have $\sum_{i=1}^s\textnormal{mult}_qB_i\leq\sum_{i=1}^s\textnormal{mult}_pB_i$. 
If it is not an equality, then we are done by induction. 
Suppose it is an equality, then we may assume that $p$ is contained in two exceptional divisors, say $E_j$ and $E_{j+1}$ for some $j$. 

\item Suppose $p$ is contained in two exceptional divisors, say $E_j$ and $E_{j+1}$ for some $j$. 
Note that 
\[a(F_1,\sF,\D) = a(E_j,\sF,\D)+a(E_{j+1},\sF,\D)-\sum_{i=1}^s(\textnormal{mult}_pB_i)b_i.\] 
Assume $a(E_j,\sF,\D)\leq a(E_{j+1},\sF,\D)$. 
If $a(F_1,\sF,\D)\geq a(E_j,\sF,\D)$, then all exceptional divisors from further blowups have discrepancies at least $a(E_j,\sF,\D)$. 
Thus, we may assume that $\sigma < a(E_j,\sF,\D)$. 
Notice that 
\[a(E_j,\sF,\D) - \sigma = \sum_{i=1}^s(\textnormal{mult}_pB_i)b_i-a(E_{j+1},\sF,\D) \geq \delta_\alpha\] 
where $\alpha = a(E_{j+1},\sF,\D)$. 
Therefore, $a(F_1,\sF,\D) \leq a(E_j,\sF,\D)-\delta_\alpha\leq 1-\delta_\alpha$. 
There are two reduced singularities $q_1$ and $q_2$ on $F_1$. 
Note that the quantities $\sum_{i=1}^s\textnormal{mult}_{q_1}B_i$ and $\sum_{i=1}^s\textnormal{mult}_{q_2}B_i$ are at most $\sum_{i=1}^s\textnormal{mult}_pB_i$.

If both are strictly inequalities, then we are done by induction. 

If one of them is an equality, say $q_1$, then the exceptional divisor $F_2$ from blowing up at $q_1$ has discrepancy 
\[a(F_2,\sF,\D)\leq a(F_1,\sF,\D) - \delta_\alpha\leq 1-2\delta_\alpha.\]
Since the minimal log discrepancy $\textnormal{mld}_x(\sF,\D)$ is positive, after finitely many blowups, the quantity $\sum_{i=1}^s\textnormal{mult}_pB_i$ will strictly decrease. 

\item Suppose $p$ is a reduced singularity contained in precisely one exceptional divisor $E_r$. 
Note that 
\begin{eqnarray*}
a(F_1,\sF,\D) &=& a(E_r,\sF,\D)-\sum_{i=1}^s(\textnormal{mult}_pB_i)b_i\\
&\leq & a(E_r,\sF,\D)-\ve\sum_{i=1}^s\textnormal{mult}_pB_i.
\end{eqnarray*} 

If $\sum_{i=1}^s\textnormal{mult}_qB_i<\sum_{i=1}^s\textnormal{mult}_pB_i$ for all $q\in F_1$, then we are done by induction. 

When $\sum_{i=1}^s\textnormal{mult}_qB_i = \sum_{i=1}^s\textnormal{mult}_pB_i$ for some $q$, if $q$ is either a smooth foliation point or contained in two exceptional divisors, then we are done. 
Otherwise, $q$ is a reduced singularity contained in precisely one exceptional divisor. 
Then the exceptional divisor $F_2$ from blowing up at $q$ has discrepancy 
\begin{eqnarray*}
a(F_2,\sF,\D) &\leq & a(F_1,\sF,\D) - \sum_{i=1}^s(\textnormal{mult}_qB_i)b_i \\
&\leq & a(E_r,\sF,\D)-2\ve\sum_{i=1}^s\textnormal{mult}_pB_i.
\end{eqnarray*}
Since the minimal log discrepancy $\textnormal{mld}_x(\sF,\D)$ is positive, after finitely many blowups, the quantity $\sum_{i=1}^s\textnormal{mult}_pB_i$ will strictly decrease. 
\end{enumerate}
\qed
\end{pf}

\begin{rmk}
When $B=\emptyset$, we have 
\[\textnormal{MLD}(2,\emptyset) = \left\{\frac{1}{n}\vert\, n\in \bN\right\}\cup\{0,-\infty\}.\]
\end{rmk}

\section{Vanishing theorem for foliations}

\begin{defn}\label{good_log_canonical}
In Theorem~\ref{main}, we say a log canonical foliation singularity is \emph{good} if it is either of the type (1)-(5), or of the type (6) and (7) such that all invariant curves are $(-1)$-$\sG$-curves and the non-invariant curve $E$ satisfying 
\[-E^2 \geq \max\{2p_a(E)-1+\deg\,[E], 2-2p_a(E)\}.\] 
\end{defn}

\begin{rmk}
All canonical foliation singularities are good log canonical foliation singularities. 
\end{rmk}

\begin{thm}\label{van_thm}
Let $f : (Y,\sG) \rw (X,\sF)$ be a proper birational morphism where $(X,\sF)$ is a foliated surface with good log canonical foliation singularities and $(Y,\sG)$ is a foliated surface with only reduced singularities. 
Then $R^if_*\cO_Y(K_\sG)=0$ for $i>0$. 
\end{thm}

\begin{pf}
We divide the proof into several steps.
\begin{enumerate}
\item We first consider $(Y,\sG)$ is the \emph{minimal resolution} of $(X,\sF)$ with exceptional divisors $E_1$, $\ldots$, $E_r$. 
Let $Z = \sum_{i=1}^ra_iE_i$ where $a_i$'s are non-negative integers. 
By the theorem on formal functions, it suffices to show that 
\[h^1(Z,\cO_Y(K_\sG)\otimes\cO_Z)=0\] 
for any effective (non-zero) divisor $Z$. 
Let $A := \sum_{i=1}^ra_i$. 
We will show the vanishing by induction on $A$. 

\item When $A = 1$, there is exactly one $i$ such that $a_i$ is positive and equals to $1$. 
Without loss of generality, we assume $a_1=1$ and $a_i=0$ for $i\geq 2$. 
\begin{claim}
$\deg(K_{E_1})-K_\sG\cdot E_1<0$. 
\end{claim}
\begin{pf}
\begin{enumerate}
\item If $E_1$ is invariant, then we have $K_\sG\cdot E_1\geq -1$ and 
\[\deg(K_{E_1})-K_\sG\cdot E_1 = -2-K_\sG\cdot E_1\leq -1<0.\]
\item If $E_1$ is non-invariant, then we have 
\begin{eqnarray*}
\deg(K_{E_1})-K_\sG\cdot E_1 &=& 2p_a(E_1)-2+E_1^2 \\
&\leq & 2p_a(E_1)-2 - (2p_a(E_1)-1+\deg\,[E_1]) \\
&=& -1-\deg\,[E_1] < 0.
\end{eqnarray*}

\end{enumerate}
\end{pf}

Thus, $h^1(E_1,\cO_Y(K_\sG)\otimes\cO_{E_1}) = h^0(E_1,\omega_{E_1}\otimes\cO_{E_1}(-K_\sG)) = 0$. 

\item 
Let $Z_\ell = Z-E_\ell$ for some $\ell$ with positive $a_\ell$, which will be determined later. 
Then we combine the following two short exact sequences 
\[\xymatrix{0 \ar[r] & \cO_Y(-Z) \ar[r] & \cO_Y(-Z_\ell) \ar[r] & \cO_{E_\ell}\otimes\cO_{Y}(-Z_\ell) \ar[r] & 0}\]
and 
\[\xymatrix{0 \ar[r] & \cO_Y(-Z_\ell)/\cO_Y(-Z) \ar[r] & \cO_Z \ar[r] & \cO_{Z_\ell}\ar[r] & 0}\]
into a short exact sequecne 
\[\xymatrix{0 \ar[r] & \cO_{E_\ell}\otimes\cO_{Y}(-Z_\ell) \ar[r] & \cO_Z \ar[r] & \cO_{Z_\ell}\ar[r] & 0}.\]
Tensoring with $\cO_Y(K_\sG)$, we get 
\[\xymatrix{0 \ar[r] & \cO_{E_\ell}\otimes \cO_Y(K_\sG-Z_\ell) \ar[r] & \cO_Y(K_\sG)\otimes\cO_Z \ar[r] & \cO_Y(K_\sG)\otimes\cO_{Z_\ell}\ar[r] & 0}.\]

Now, by induction hypothesis, we have $h^1(Z_\ell,\cO_Y(K_\sG)\otimes\cO_{Z_\ell})=0$. 
Therefore, it suffices to show $h^1(E_\ell,\cO_{E_\ell}\otimes \cO_Y(K_\sG-Z_\ell))=0$. 
This follows from 
\[K_\sG\cdot E_\ell-Z_i\cdot E_\ell>\deg(K_{E_\ell}) = K_Y\cdot E_\ell+E_\ell^2,\]
which is equivalent to  
\[(K_\sG-K_Y-Z)\cdot E_\ell>0.\]

\item Let $K_\sG-K_Y\equiv_f\sum_{i=1}^rb(E_i)E_i$ for some $b(E_i)\in\bR$. 
We will denote $b(E_i)$ by $b_i$ as well. 
\begin{claim}
$b(E_i)\leq 1$ for all $i$. 
Moreover, either $b(E_i)=1$ for all $i$ or $b(E_i)<1$ for all $i$. 
\end{claim}
\begin{pf}
\begin{enumerate}
\item If $E_j$ is invariant, then we have 
\begin{eqnarray*}
(K_\sG-K_Y)\cdot E_j &=& \textnormal{Z}(\sG,E_j)-2+E_j^2+2 \\
&=&\textnormal{Z}(\sG,E_j)+E_j^2 \geq \left(\sum_{i=1}^rE_i\right)\cdot E_j.
\end{eqnarray*}
\item If $E_j$ is non-invariant, then we have 
\begin{eqnarray*}
(K_\sG-K_Y)\cdot E_j &=& -E_j^2 + E_j^2+2-2p_a(E_j) \\
&> & \deg\,[E_j]+E_j^2 = \left(\sum_{i=1}^rE_i\right)\cdot E_j.\end{eqnarray*}
\end{enumerate} 
By Lemma~\ref{key_lem}, we have that $b_i\leq 1$ for all $i$. 
Moreover, $\bigcup_{i=1}^rE_i$ is connected, we have either $b_i=1$ for all $i$ or $b_i<1$ for all $i$. 
\end{pf}

\begin{claim}
We have $b(E_i)\geq 0$ for all $i$ except for the non-invariant $E_0$ of (6) and (7) in Theorem~\ref{main}. 
\end{claim}
\begin{pf}
Let $2D$ be the sum of all $(-1)$-$\sG$-curves whose self-intersections are $-2$. 
Then we have $(K_\sG-K_Y)\cdot E_i\leq D\cdot E_i$ for all $i$. 
Thus, by Lemma~\ref{key_lem}, we have that $b_i\geq 0$ for all $i$. 

In the case of type (6) and (7), we note that 
\[(K_\sG-K_Y)\cdot E_0\leq 2-2p_a(E_0) \leq -E_0^2\]
and 
\[(K_\sG-K_Y)\cdot E_j = 1+E_j^2\leq -1 = -E_0\cdot E_j\] 
for any invariant curve $E_j$. 
Then, by Lemma~\ref{key_lem}, we have $b(E_j)\geq 0$ for all invariant curves $E_j$. 
\end{pf}

\item\label{anti-ample} 
Let $W = \sum_{i=1}^rc_iE_i$ be an $f$-anti-ample effective divisor where $c_i>0$ for all $i$. 
We have following two cases:
\begin{enumerate}
\item If $a_i-b_i$ is positive for some $i$, then we define $\alpha\in\bR_{\geq 0}$ to be the maximal number such that $\alpha(a_i-b_i)\leq c_i$ for all $i$. 
Then there is in index $\ell$ such that $c_\ell-\alpha(a_\ell-b_\ell)=0$, and therefore 
\[\big(W+\alpha(K_\sG-K_Y-Z)\big)\cdot E_\ell = \left(\sum_{i=1}^r \big(c_i-\alpha(a_i-b_i)\big)E_i\right)\cdot E_\ell\geq 0.\]
Hence, we have 
\[(K_\sG-K_Y-Z)\cdot E_\ell\geq\frac{-1}{\alpha}W\cdot E_\ell>0.\]

\item If $a_i-b_i\leq 0$ for all $i$, then $a_i\leq b_i\leq 1$ for all $i$. 
Thus, $b_i=1$ for all $i$, otherwise, we have $b_i<1$ for all $i$, and then $a_i=0$ for all $i$ since $a_i$'s are non-negative integers. 
\end{enumerate}

\item Suppose $\bigcup_{i=1}^rE_i$ forms a dual graph for a canonical foliation singularity. 
Then we have $b(E_i)\geq 0$ for all $i$. 
If $b_i<1$ for all $i$, then, from $c_\ell-\alpha(a_\ell-b_\ell)=0$, we have $a_\ell = b_\ell+\frac{c_\ell}{\alpha}>0$. 

If $b_i=1$ for all $i$, then $a_i\leq 1$ for all $i$ and $\bigcup_{i=1}^rE_i$ is an elliptic Gorenstein leaf $\Gamma$. 
By \cite[Fact III.0.4 and Theorem IV.2.2]{mcquillan2008canonical}, we have $\cO_\G(K_\sG)$ is not torsion and has degree $0$. 
Also $\G$ is Cohen-Macaulay with trivial dualizing sheaf, then by Serre duality, we have 
\[h^1(\G,\cO_Y(K_\sG)\otimes\cO_\G) = h^0(\G,\cO_Y(-K_\sG)\otimes\cO_\G) = 0.\]
Thus, we have shown the case when $a_i=1$ for all $i$. 
Therefore, when $a_i=0$ for some $i$, we have an index $j$ such that 
\[(K_\sG-K_Y-Z)\cdot E_j = \left(\sum_i(b_i-a_i)E_i\right)\cdot E_j> 0.\]

\item Now suppose $\bigcup_{i=1}^rE_i$ forms a dual graph for a non-canonical foliation singularity. 
Note that we have $b_i<1$ for all $i$. 
If $E_\ell$ in step~(\ref{anti-ample}) is invariant, then from $c_\ell-\alpha(a_\ell-b_\ell)=0$, we have $a_\ell = b_\ell+\frac{c_\ell}{\alpha}>0$. 
Thus, it suffices to consider the case when $E_\ell$ is non-invariant and $a_\ell = 0$. 
Then the support of $Z = \sum_iZ_i = \sum_ia_iE_i$ is the disjoint union of $(-1)$-$\sG$-curves $E_i$ where $Z_i = a_iE_i$. 
Note that 
\[H^1(Z,\cO_Y(K_\sG)\otimes\cO_Z)=\bigoplus_i H^1(Z_i,\cO_Y(K_\sG)\otimes\cO_{Z_i}).\]
Since the contraction of $E_i$ introduces a terminal foliation singularity, we have $h^1(Z_i,\cO_Y(K_\sG)\otimes\cO_{Z_i}) = 0$ for all $i$, and thus $h^1(Z,\cO_Y(K_\sG)\otimes\cO_Z)=0$. 

\item Now we consider the general resolution. 
Note that any general resolution $f : (Y,\sG) \rw (X,\sF)$ factors through the minimal resolution $h : (Z,\sH) \rw (X,\sF)$. 
Let $f = h\circ g$. 
So $g$ is the composition of blowups. 

We have shown that $R^1h_*\cO_Z(K_\sH) = 0$. 
Notice that $g_*\cO_Y(K_\sG) = \cO_Z(K_\sH)$. 
By Grothendieck spectral sequence, we have the following exact sequence:
\[\xymatrix{0 \ar[r] & R^1h_*(g_*\cO_Y(K_\sG)) \ar[r] & R^1f_*\cO_Y(K_\sG) \ar[r] & h_*R^1g_*\cO_Y(K_\sG).}\]
Since $R^1h_*(g_*\cO_Y(K_\sG)) = R^1h_*\cO_Z(K_\sH)=0$ and $R^1g_*\cO_Y(K_\sG)=0$ by Lemma~\ref{blowup_van_thm}, we have $R^1f_*\cO_Y(K_\sG)=0$. 
\qed
\end{enumerate}
\end{pf}

\begin{lem}\label{blowup_van_thm}
Let $(X,\sF,p)$ is a germ of a foliated surface with at worst reduced singularity at $p$. 
Let $\pi : (Y,E) \rw (X,p)$ be a blowup at $p$ and $\sG$ be the pullback foliation of $\sF$. 
Then we have $R^1\pi_*\cO_Y(K_\sG) = 0$.
\end{lem}
\begin{pf}
Note that $K_\sG = \pi^*K_\sF + a(E)E$ where $a(E)$ is either $0$ or $1$. 
Also, since $p$ is a smooth point of $X$, we have $R^1\pi_*\cO_Y=0$. 

If $a(E)=0$, then, by the projection formula, we have 
\[R^1\pi_*\cO_Y(K_\sG) = R^1\pi_*\pi^*\cO_X(K_\sF) = R^1\pi_*\cO_Y\otimes\cO_X(K_\sF)=0.\]

If $a(E)=1$, then we consider the following short exact sequence:
\[\xymatrix{0 \ar[r] & \pi^*\cO_X(K_\sF) \ar[r] & \cO_Y(K_\sG) \ar[r] & \cO_E(E) \ar[r] & 0}\]
Pushing forward via $\pi$, we obtain the exact sequence 
\[\xymatrix{0=R^1\pi_*\pi^*\cO_X(K_\sF) \ar[r] & R^1\pi_*\cO_Y(K_\sG) \ar[r] & R^1\pi_*\cO_E(E) \ar[r] & 0}.\]
Since $R^1\pi_*\cO_E(E) = H^1(\bP^1,\cO_{\bP^1}(-1)) = 0$, we have $R^1\pi_*\cO_Y(K_\sG) = 0$.
\qed
\end{pf}

\bibliographystyle{amsalpha}
\addcontentsline{toc}{chapter}{\bibname}
\normalem
\bibliography{LC_Foliation}

\end{document}